\documentclass[10pt,draftcls,onecolumn]{IEEEtran}  

\IEEEoverridecommandlockouts

\usepackage{verbatim}
\usepackage{cite}
\usepackage{multicol}
\usepackage{amsmath}
\usepackage{amssymb}
\usepackage{graphicx,psfrag}
\usepackage{caption}
\usepackage{subcaption}
\usepackage{tikz}

\usepackage{float}

\title{Local stability and Hopf bifurcation analysis for Compound TCP}

\author{Debayani Ghosh, Krishna Jagannathan and Gaurav Raina
\thanks{D.~Ghosh, K.~Jagannathan and G.~Raina are with the Department of Electrical Engineering, IIT
Madras, Chennai 600036, India. Email: \{ee12s052,   krishnaj, gaurav\}@ee.iitm.ac.in
        }
}

\begin{document}

\maketitle
\begin{abstract}
We conduct a local stability and Hopf bifurcation analysis for Compound TCP, with small Drop-tail buffers, in three topologies. The first topology consists of two sets of TCP flows having different round trip times, and feeding into a core router. The second topology corresponds to two queues in tandem, and consists of two distinct sets of TCP flows, regulated by a single edge router and feeding into a core router. The third topology comprises of two distinct sets of TCP flows, regulated by two separate edge routers, and feeding into a common core router. For each of these cases, we conduct a detailed local stability analysis and obtain conditions on the network and protocol parameters to ensure stability. If these conditions get marginally violated, our analysis shows that the underlying systems would lose local stability via a Hopf bifurcation. After exhibiting a Hopf, a key concern is to determine the asymptotic orbital stability of the bifurcating limit cycles. We present a detailed analytical framework to address the stability of the limit cycles, and the type of the Hopf bifurcation by invoking Poincar\'{e} normal forms and the center manifold theory. We conduct packet-level simulations to highlight the existence and stability of the limit cycles in the queue size dynamics. 
\end{abstract}

\section{Introduction}
Network performance, and end-to-end latency are affected by a combination of the choice of TCP, the size of router buffers, and the choice of queue management implemented in Internet routers \cite{Cerf},~\cite{Gettys},~\cite{Nichols}. A major portion of Internet traffic is controlled by the Transmission Control Protocol (TCP)~\cite{Ha},~\cite{Padhye}. There have been proposals for different flavours of TCP and queue management strategies. However, Compound TCP~\cite{Tan} is the default protocol in Windows, and a simple Drop-Tail queue management is commonly implemented in Internet routers. It has been argued that the default large buffer dimensioning rule for router buffers, combined with Drop-Tail, leads to excessive delays in the Internet \cite{Gettys}.

\begin{figure}
\begin{center}
   \psfrag{0}{\begin{scriptsize}$0$\end{scriptsize}}
  \psfrag{15}{\begin{scriptsize}$15$\end{scriptsize}}
  \psfrag{100}{\begin{scriptsize}$100$\end{scriptsize}}
  \psfrag{125}{\begin{scriptsize}$125$\end{scriptsize}}
      \psfrag{aaaa}{\begin{scriptsize}Buffer size = 15 pkts\end{scriptsize}}
  \psfrag{bbbb}{\begin{scriptsize}Buffer size = 100 pkts\end{scriptsize}}
\psfrag{cccc}{\hspace{-4mm} Queue size (pkts)}
\psfrag{xyz}{\hspace{-8mm}Time (seconds)}
  \psfrag{mmmm}{\begin{scriptsize}Round trip time = $10$ ms\end{scriptsize}}
   \psfrag{nnnn}{\begin{scriptsize}\hspace{-1mm}Round trip time = $200$ ms\end{scriptsize}}

  \includegraphics[width=2.5in,height=3.5in,angle=270]{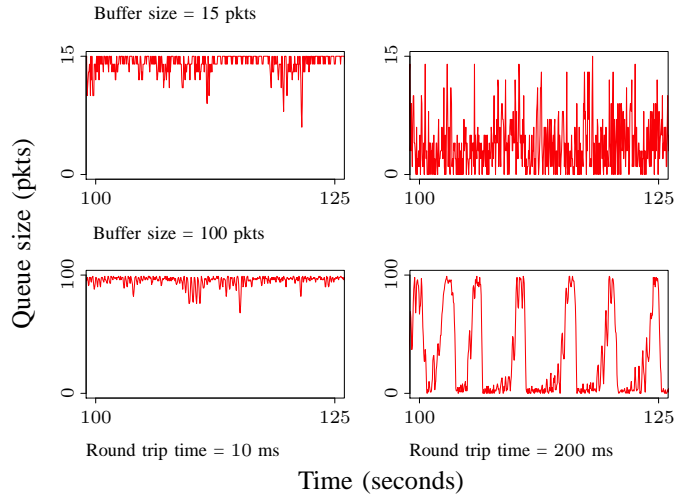}
  \caption{ \emph{Long-lived flows}. 60 long-lived Compound flows over a 2 Mbps link, and feeding into a core router with link capacity 100 Mbps. Observe the emergence of limit cycles in the queue at the core router, for larger buffer thresholds, and larger round trip times.}\vspace{-8mm}
  \label{ns_single}
  \end{center}
\end{figure}

In our recent work~\cite{Jagannathan}, we conducted a performance evaluation of Compound TCP, in a small buffer regime, with particular emphasis on buffer thresholds. One of the key insights obtained therein was the two-fold advantage of having small router buffers. In particular, our analysis showed that small buffers are favourable for ensuring the stability of the system, in addition to reducing queuing delays. Furthermore, our analysis identified that the underlying dynamical systems undergo a Hopf bifurcation, and transit from a locally stable into an unstable regime as the buffer size increases. The Hopf bifurcation alerts us to the emergence of isolated periodic orbits, termed as limit cycles, as a parameter crosses a certain critical value. In addition, we repeatedly observed limit cycles in the queue size dynamics, in numerous packet-level simulations. Fig.~\ref{ns_single} portrays one such instance; indeed, it captures the emergence of limit cycles in the queue size of the core router in a single bottleneck topology. This motivates us to develop an analytical framework under which the emergence of these non-linear oscillations can be better understood. To that end, in this paper, we provide a complete analytical characterisation of the type of the Hopf bifurcation, and prove the orbital stability of the emergent limit cycles.           

We consider three different topologies, and focus on analysing the dynamical properties of a fluid model of Compound TCP in conjunction with small Drop-Tail buffers. Our fluid model takes the form of a non-linear, time-delayed dynamical system. The first topology is a generalisation of the single bottleneck topology studied in~\cite{Jagannathan}, and consists of two sets of TCP flows having different round trip times, and feeding into a core router (see Fig.~\ref{Figure.1}(a)). The second topology corresponds to two queues in tandem, and consists of two distinct sets of TCP flows, regulated by a single edge router and feeding into a core router (see Fig.~\ref{Figure.1}(b)). The third topology comprises of two distinct sets of TCP flows, regulated by two separate edge routers, and feeding into a common core router (see Fig.~\ref{Figure.1}(c)). 

For each of these cases, we first conduct a local stability analysis and outline necessary and sufficient conditions for local stability, with two simplifying assumptions. In the first scenario, we assume that the network parameters are the same, and that both sets of Compound TCP flows have equal round trip times. In the second scenario, we assume the network parameters to be heterogeneous, and the round trip time of one set of TCP flows to be much larger as compared to the other. If the local stability conditions get marginally violated, our analysis shows that the underlying systems would lose local stability via a Hopf bifurcation. Motivated by this insight, we then analyse only the third topology in greater detail, to better understand the impact of heterogeneous system parameters on local stability. We numerically show through DDE-BIFTOOL~\cite{DDE1},~\cite{DDE2} that, even in the presence of heterogeneous network parameters and different round trip times, the dynamical system undergoes a Hopf bifurcation which leads to the emergence of limit cycles. 

As argued in~\cite{Jagannathan}, the emergence of limit cycles in the system dynamics could have a number of detrimental consequences -- for example it could lead to the synchronisation of TCP windows, result in a loss in link utilisation, and cause the downstream traffic to be bursty. Hence, it becomes imperative to study these limit cycles in further detail. To that end, an important contribution of this paper lies in providing an analytical framework to determine the asymptotic orbital stability of the emerging limit cycles. Using Poincar\'{e} normal forms and the center manifold theory, we show that the Hopf bifurcation is indeed supercritical, and hence leads to the emergence of orbitally stable limit cycles. To corroborate our analytical insights, we conduct some packet-level simulations in NS2~\cite{ns2},  to highlight the existence and stability of the limit cycles in the queue size dynamics. Notably, instead of treating any particular system parameter as the bifurcation parameter, we choose a suitably motivated exogenous, non-dimensional parameter as the bifurcation parameter to aid our analysis. The two main advantages of this are: first, it enables us to capture the effects of different system parameters on the system stability in a unified manner and secondly, we need not be concerned about the dimension of the bifurcation parameter.                  

The rest of the paper is as organised as follows. In section \ref{Models}, we outline the governing fluid models for the three cases we consider. Section \ref{Localstability} deals with local stability analysis of the fluid models. In Section \ref{sec:hopf}, we provide an analytical framework to determine the asymptotic orbital stability of the bifurcating limit cycles, and to characterise the type of the Hopf bifurcation. Packet-level simulations are presented in Section~\ref{simulations} to corroborate some of the analytical insights. Finally, in section \ref{conclusions}, we summarise our contributions.      
\section{Models}
\label{Models}

In this section, we consider two distinct sets of TCP flows having different round trip times $\tau_1$ and $\tau_2$ in three topologies. For our analysis of these models, we primarily focus on long-lived flows. We assume that both sets of TCP flows can be of different flavours and hence, can have different increase and decrease rules to govern the evolution of the corresponding window sizes. Let the average window sizes of the two sets of flows be $w_{1}(t)$ and $w_{2}(t)$  respectively. For each acknowledgement received, the average window sizes increase by $i_1(w_1(t))$ and $i_2(w_2(t))$, and for each packet loss detected, the average window sizes decrease by  $d_1(w_1(t))$ and $d_2(w_2(t))$ respectively. Note that, the increase and decrease functions for a particular TCP flavour depend on the protocol parameters. Further, the loss probability at the routers is governed by the corresponding AQM strategy. 

\subsection{Fluid models for TCP}
Now, we briefly outline the fluid models for the evolution of the average window sizes of the two sets of TCP flows in the congestion avoidance phase for three topologies. 
\subsection*{Case I}
This model consists of a single bottleneck link with two distinct sets of TCP flows feeding into a  common core router, as shown in Fig.~\ref{Figure.1}(a). The core router has a buffer size of $B$, with link capacity $C$. Thus, for generalised TCP flows, the non-linear, time-delayed, fluid model of the system is given by the following equations:
\\
\begin{align}
\frac{dw_{j}(t)}{dt} = \frac{w_{j}(t-\tau_{j})}{\tau_{j}}\bigg(i_j\left(w_{j}(t)\right)\Big(1-q(t,\tau_1,\tau_2)\Big) - d_j\left((w_{j}(t)\right)q(t,\tau_1,\tau_2)\bigg), \hspace{2ex}j=1,2,
\label{eq:modelb_1}
\end{align}
where $q(t,\tau_1,\tau_2)$ represents the packet loss probability at the core router, and depend on the sending rates of both sets of TCP flows.
\subsection*{Case II}
This model consists of two distinct sets of TCP flows, regulated by a single edge router and feeding into a common core router, as shown in Fig.~\ref{Figure.1}(b). The buffer sizes of the core router and the edge routers are $B_1$ and $B_2$, with link capacities $C_1$ and $C_2$ respectively. Thus, for generalised TCP flows, the non-linear, time-delayed, fluid model of the system is given by the following differential equations:

\begin{align}
\frac{dw_{j}(t)}{dt} = \frac{w_{j}(t-\tau_{j})}{\tau_{j}}\bigg(i_j\left(w_{j}(t)\right)\Big(1-q_1(t,\tau_1,\tau_2)-q_2(t,\tau_1,\tau_2))\Big) - d_j\left((w_{j}(t)\right)\left(q_1(t,\tau_1,\tau_2)+q_2(t,\tau_1,\tau_2)\right)\bigg),
\label{eq:modelb_2}
\end{align}
for $j=1,2,$ and $q_1(t,\tau_1,\tau_2)$ and $q_2(t,\tau_1,\tau_2)$ denote the packet loss probabilities at the edge router and the core router respectively.
\subsection*{Case III}
\begin{figure*}
\vspace{2mm}
{
\begin{tikzpicture}[scale=0.5]
\draw (0,0) -- (-1.5,0.75);
\draw (0,0) -- (-1.5,0.25);
\draw (0,0) -- (-1.5,-0.25);
\draw (-1.5,-0.75) -- (0,0);
\draw (0,-0.75) -- (3.85,-0.75);
\draw (3.85,-0.75) -- (3.85,0.75);
\draw (3.85,0.75) -- (0,0.75);
\draw (0,0.75) -- (0,-0.75);
\draw (3.65,-0.75) -- (3.65,0.75);
\draw (3.4,-0.75) -- (3.4,0.75);
\draw (3.15,-0.75) -- (3.15,0.75);
\draw (4.65,0) circle(0.75 cm);
\draw [loosely dotted,very thick] (3.0,0) -- (2.4,0);
\draw[->] (5.4,0) -- (6.2,0);
\draw (5.8,0) -- (5.8,2.2);
\draw (5.8,2.2) -- (-1.75,2.2);
\draw (-1.75,2.2) -- (-1.75,0.5);
\draw [->](-1.75,0.5) -- (-1.5,0.5);
\draw (5.8,0) -- (5.8,-2.2);
\draw (5.8,-2.2) -- (-1.75,-2.2);
\draw (-1.75,-2.2) -- (-1.75,-0.5);
\draw [->](-1.75,-0.5) -- (-1.5,-0.5);
\node at (1,0) {$B$};
\node at (4.65,0) {$C$};
\node at (1.5, 1.75) {$\tau_1$};
\node at (1.5,-1.9) {$\tau_2$};
\node at (1.75,-3) {(a)};

\draw (8.5,0) -- (8,0.75);
\draw (8.5,0) -- (8,0.25);
\draw (8.5,0) -- (8,-0.25);
\draw (8.5,0) -- (8,-0.75);
\draw (8.5,-0.75) -- (10.7,-0.75);
\draw (10.7,-0.75) -- (10.7,0.75);
\draw (10.7,0.75) -- (8.5,0.75);
\draw (8.5,0.75) -- (8.5,-0.75);
\draw (10.55,-0.75) -- (10.55,0.75);
\draw (10.35,-0.75) -- (10.35,0.75);
\draw [loosely dotted,very thick,] (10.25,0) -- (9.75,0);
\draw (12.25,0) -- (13,0);
\draw (13,-0.75) -- (15.2,-0.75);
\draw (15.2,-0.75) -- (15.2,0.75);
\draw (15.2,0.75) -- (13,0.75);
\draw (13,0.75) -- (13,-0.75);
\draw (15.05,-0.75) -- (15.05,0.75);
\draw (14.85,-0.75) -- (14.85,0.75);
\draw (11.5,0) circle(0.75 cm);
\draw (16,0) circle(0.75 cm);
\draw [loosely dotted,very thick] (14.8,0) -- (14.4,0);
\draw[->] (16.75,0) -- (17.5,0);
\draw (17.125,0) -- (17.125,2.2);
\draw (17.125,2.2) -- (7.75,2.2);
\draw (7.75,2.2) -- (7.75,0.5);
\draw [->](7.75,0.5) -- (8,0.5);
\draw (17.125,0) -- (17.125,-2.2);
\draw (17.125,-2.2) -- (7.75,-2.2);
\draw (7.75,-2.2) -- (7.75,-0.5);
\draw [->](7.75,-0.5) -- (8,-0.5);
\node at (9.2,0) {$B_1$};
\node at (11.5,0) {$C_1$};
\node at (13.7,0) {$B_2$};
\node at (16.0,0) {$C_2$};
\node at (12.5, 1.75) {$\tau_1$};
\node at (12.5, -1.9) {$\tau_2$};
\node at (12.5,-3){(b)};


\draw (8.5,0) -- (8,0.75);
\draw (8.5,0) -- (8,0.25);
\draw (8.5,0) -- (8,-0.25);
\draw (8.5,0) -- (8,-0.75);
\draw (19.8,0.5) -- (22.5,0.5);
\draw (22.5,0.5) -- (22.5,1.7);
\draw (22.5,1.7) -- (19.8,1.7);
\draw (19.8,1.7) -- (19.8,0.5);
\draw (22.25,0.5) -- (22.25,1.7);
\draw (22,0.5) --(22,1.7);
\draw [loosely dotted,very thick] (22,1.1) -- (21.5,1.1);
\node at (20.4,1.1) {$B_1$};
\node at (23.15,1.1) {$C_1$};

\draw (19.8,-0.5) -- (22.5,-0.5);
\draw (22.5,-0.5)-- (22.5,-1.7);
\draw (22.5,-1.7) -- (19.8,-1.7);
\draw (19.8,-1.7) -- (19.8,-0.5);
\draw (22.25,-0.5) -- (22.25,-1.7);
\draw (22,-0.5) --(22,-1.7);
\draw [loosely dotted,very thick] (22,-1.1) -- (21.5,-1.1);
\draw [loosely dotted,very thick] (27.6,0) -- (27,0);
\node at (20.4,-1.1) {$B_2$};
\node at (23.15,-1.1) {$C_2$};
\node at (24.5,1.75) {$\tau_1$};
\node at (24.5,-1.9) {$\tau_2$};	
\node at (26,0) {$B$};
\node at (28.85,0) {$C$};
\node at (24.5,-3) {(c)};

\draw (23.15,1.1) circle(0.6 cm);
\draw (23.15,-1.1) circle(0.6 cm);
\draw (23.75,1.1) -- (25.5,0);
\draw (23.75,-1.1) -- (25.5, 0);

\draw (25.5,0.6) -- (28.2,0.6);
\draw (28.2,0.6) -- (28.2,-0.6);
\draw (28.2,-0.6) -- (25.5,-0.6);
\draw (25.5,-0.6) -- (25.5,0.6);
\draw (27.95,-0.6) -- (27.95,0.6);
\draw (27.7,-0.6) -- (27.7,0.6);
\draw (27.7,-0.6) -- (27.7,0.6);

\draw (28.85,0) circle(0.6 cm);
\draw[->] (29.45,0) -- (31,0);
\draw (18.8,1.7) -- (19.8,1.1);
\draw (18.8,0.5) -- (19.8,1.1);
\draw (18.8,-0.5) -- (19.8,-1.1);
\draw (18.8,-1.7) -- (19.8,-1.1);
\draw (30.225,0) -- (30.225,2.2);
\draw (30.225,2.2) -- (18.5,2.2);
\draw (18.5,2.2) -- (18.5,1.1);
\draw[->] (18.5,1.1) -- (18.9,1.1);
\draw (30.225,0) -- (30.225,-2.2);
\draw (30.225,-2.2) -- (18.5,-2.2);
\draw (18.5,-2.2) -- (18.5,-1.1);
\draw[->] (18.5,-1.1) -- (18.9,-1.1);

\end{tikzpicture}

}
\caption{\emph{Schematic diagrams of three topologies.} (a) Case I, a single bottleneck topology (b) Case II, two routers in tandem and (c) Case III, two routers feeding into one core router.}
\label{Figure.1}
\end{figure*}
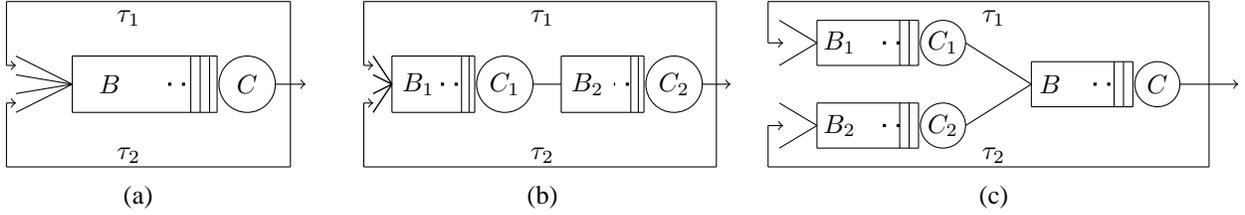

This model consists of  two distinct sets of TCP flows, regulated by two edge routers and feeding into a common core router, as shown in Fig.~\ref{Figure.1}(c). The buffer size at the core router is $B$, with link capacity $C$. The buffer sizes for the edge routers are $B_1$ and $B_2$, with link capacities $C_1$ and $C_2$ respectively. Thus, for generalised TCP flows, the non-linear, time-delayed, fluid model of the system is given by the following equations:\\
\begin{align}
\frac{dw_{j}(t)}{dt} = \frac{w_{j}(t-\tau_{j})}{\tau_{j}}\bigg(i_j\left(w_{j}(t)\right)\Big(1-p_{j}(t-\tau_{j})-q(t,\tau_{1},\tau_{2})\Big) - d_j\left((w_{j}(t)\right)\Big(p_{j}(t-\tau_{j})+q(t,\tau_{1},\tau_{2})\Big)\bigg),
\label{eq:modelb_3}
\end{align}
for $j=1,2.$ The loss probabilities at the two edge routers are $p_1(t)$ and $p_2(t)$. The loss probability at the core router is denoted by $q(t, \tau_1, \tau_2)$. Recall that, the increase and decrease functions are specific to the choice of a particular flavour of TCP. Specifically,~\cite{Raja} has summarised the increase and decrease functions for different TCP flavours including Compound. Since our primary focus is on Compound TCP, we state the increase and decrease functions for Compound as follows:
\begin{equation}
\label{eq:increasedecrease}
i(w(t))=\alpha\left(w(t)\right)^{k-1}, \hspace{1ex}\text{and}\hspace{1ex} d(w(t))=\beta w(t).
\end{equation} 
Here, $\alpha$, $k$ are the increase parameters and $\beta$ is the decrease parameter. The default values of these parameters are $\alpha=0.125$, $k=0.75$ and $\beta=0.5$~\cite{Tan}.
\subsection{Packet loss probability}
In this paper, we mainly focus on small buffers with Drop-Tail queue policy for the local stability analysis of the non-linear fluid models of TCP given by \eqref{eq:modelb_1}, \eqref{eq:modelb_2} and \eqref{eq:modelb_3}. We first consider the scenario where a \emph{large} number of long-lived TCP flows having a common round trip time of $\tau$ feed into a router having a buffer size of $B$. The bottleneck link has a capacity $C$. In this scenario, we can approximate the packet loss probability of the router by the blocking probability of an $M/M/1/B$ queue~\cite{Raja}. This gives rise to the following fluid model:
\begin{align}
\label{eq:loss}
p(t)=\bigg(\frac{w(t)}{C\tau}\bigg)^B,
\end{align} 
where $w(t)$ represents the average window size of the TCP flows. Using \eqref{eq:loss}, we can then obtain the functional forms of packet loss probabilities for the three scenarios, which we briefly outline as follows:
\subsubsection*{Case I}The fluid model for the loss probability at the core router is given by
\begin{align}
q(t)=\Bigg(\frac{w_1(t) / \tau_1 + w_2(t) / \tau_2}{C}\Bigg)^B.
\label{eq:case1_loss}
\end{align} 
\subsubsection*{Case II} The fluid models for the loss probabilities are:
\begin{align}
q_1(t)= \left(\frac{w_{1}(t) / \tau_{1}+w_{2}(t)/ \tau_{2}}{C_2}\right)^{B_2}  {\rm and}, \hspace{1ex} q_2(t)= \left(\frac{w_{1}(t) / \tau_{1}+ w_{2}(t) /\tau_{2}}{C_1}\right)^{B_1}.
\label{eq:case2_loss}
\end{align}
\subsubsection*{Case III} Using \eqref{eq:loss}, we can approximate the loss probabilities at various routers as below:
\begin{align}
 p_{1}(t) &= \left(\frac{w_{1}(t)}{C_{1}\tau_{1}}\right)^{B_{1}} \hspace{1ex},\hspace{1ex} p_{2}(t)=\left(\frac{w_{2}(t)}{C_{2}\tau_{2}}\right)^{B_{2}},\ {\rm and} \notag \\
q(t) &= \left(\frac{w_{1}(t)/\tau_{1}+w_{2}(t)/\tau_{2}}{C}\right)^{B}.
\label{eq:case3_loss}
\end{align}
Using these functional forms, we now proceed to perform a local stability and bifurcation  analysis for the systems given by \eqref{eq:modelb_1}, \eqref{eq:modelb_2} and \eqref{eq:modelb_3}. This would enable us to understand the dynamical properties of the coupled system of Compound TCP with Drop-Tail queue policy to a greater detail.  
\section{Local Stability Analysis}
\label{Localstability}
Note that, to perform a local stability and bifurcation analysis for the non-linear models \eqref{eq:modelb_1}, \eqref{eq:modelb_2} and  \eqref{eq:modelb_3}, we need to choose an appropriate bifurcation parameter. It can be easily seen that both protocol and network parameters affect the stability of the systems. To that end, instead of treating any of the system parameters as a bifurcation parameter, we introduce an exogenous non-dimensional parameter $\kappa$ as the bifurcation parameter. We choose the non-dimensional parameter in such a manner that it does not affect the equilibrium of the system. Recall that, to conduct the local stability analysis, we primarily focus on Compound TCP with Drop-Tail queues in the small buffer regime. For mathematical tractability, we assume that both sets of TCP flows in all three topologies are regulated by Compound with identical protocol parameters. Further, we consider two simplifying assumptions as briefly outlined below:
\\\\
\emph{Scenario 1:} All network parameters are the same, \emph{i.e.}, $B_1=B_2=B$, and $C_1=C_2=C$. Further, the round trip times of both TCP flow sets are identical,  \emph{i.e.}, $\tau_1=\tau_2=\tau.$
\\\\
\emph{Scenario 2:} In this scenario, we assume that all network parameters are distinct and the round trip time of one set of TCP flows is negligible and much smaller as compared to the round trip time of the other set, \emph{i.e}, $\tau_1>>\tau_2$ and $\tau_2\approx0$. Under this assumption, the dynamics of the second set of TCP flows appear almost instantaneous.\\

We now proceed to conduct a detailed local stability analysis to obtain bounds on network, and protocol parameters to ensure stability, for the systems given by \eqref{eq:modelb_1}, \eqref{eq:modelb_2} and \eqref{eq:modelb_3}.
\subsection*{Case I}
The schematic diagram of the topology is presented in Fig.~\ref{Figure.1}(a).
\subsubsection*{Scenario 1} With this assumption, the first model reduces to a single bottleneck link with only one set of TCP flows having round trip time $\tau_1=\tau_2=\tau$. Hence, with the non-dimensional bifurcation parameter $\kappa$, system \eqref{eq:modelb_1} reduces to the following non-linear, first-order, time-delayed differential equation:
\begin{align}
\frac{dw(t)}{dt} =  \kappa\frac{w(t-\tau)}{\tau}\bigg(i\left(w(t)\right)\Big(1-q(w(t-\tau))\Big) - d\left((w(t)\right)q(w(t-\tau))\bigg),
\label{eq:modelb_1reduced1}
\end{align}
where $w(t)$ is the \emph{average} window size of the TCP flows. The non-trivial equilibrium $w^{*}$ of system \eqref{eq:modelb_1reduced1} satisfies the following equation
\begin{equation}
i(w^{*})=d(w^{*})q(w^{*}).
\end{equation}
Note that, under the first assumption, the fluid model for the packet loss probability at the core router, given by \eqref{eq:case1_loss} reduces to
\begin{align}
q(w^{*})=\left(\frac{w^{*}}{C\tau}\right)^B,
\end{align} 
at equilibrium. A necessary and sufficient condition for this model, with Compound TCP in the small buffer regime is~\cite{Raja}
\begin{align}
\kappa\alpha \left(w^{*}\right)^{k-1}\sqrt{B^{2}-\left(\left(k-2\right)\left(1-q(w^{*})\right)\right)^{2}} < \cos^{-1}\left(\frac{\left(k-2\right)\left(1-q(w^*)\right)}{B}\right).
\label{eq:condition_case1}
\end{align}
\subsubsection*{Scenario 2}
With the introduction of the non-dimensional parameter $\kappa$, system \eqref{eq:modelb_2} becomes
\begin{align}
\dot{w}_1(t) & = \kappa\frac{w_{1}(t-\tau_{1})}{\tau_{1}}\bigg(i\left(w_{1}(t)\right)\left(1-q(t,\tau_1,\tau_2)\right) - d\left((w_{1}(t)\right)q(t,\tau_{1},\tau_{2})\bigg),\notag \\
\dot{w}_2(t) & =  \kappa\frac{w_{2}(t)}{\tau_{2}}\bigg(i\left(w_{2}(t)\right)\left(1-q(t,\tau_1,\tau_2)\right) - d\left((w_{2}(t)\right)q(t,\tau_{1},\tau_{2})\bigg).
 \label{eq:modelb_1reduced2}
\end{align}
Suppose $(w_1^{*}, w_{2}^{*})$ is a non-trivial equilibrium of \eqref{eq:modelb_1reduced2} and let $u_{1}(t)=w_{1}(t)-w_{1}^{*}$ and $u_{2}(t)=w_{2}(t)-w_{2}^{*}$ be small perturbations about $w_{1}^{*}$ and $w_{2}^{*}$ respectively. Linearising \eqref{eq:modelb_1reduced2} about this equilibrium, we obtain
\begin{align}
\label{eq:linearmodelb_1reduced2}
&\dot{u}_1(t) = -\kappa\left(\mathcal{M}_{1}u_{1}(t)+\mathcal{N}_{1}u_{1}(t-\tau_{1})+\mathcal{P}_{1}u_{2}(t)\right),\notag\\
&\dot{u}_2(t) = -\kappa\left(\big(\mathcal{M}_{2}+\mathcal{N}_{2}\big)u_{2}(t)+\mathcal{P}_{2}u_{1}(t-\tau_{1})\right),
\end{align}
where, the increase and decrease functions for Compound TCP given by \eqref{eq:increasedecrease}, and the functional form of the loss probability at the core router given by~\eqref{eq:case1_loss} yield the following coefficients: 
\begin{align}
\label{eq:coefficients}
\mathcal{M}_{j} &=-\frac{\alpha}{\tau_{j}}\left(k-2\right)\ \left(w_{j}^{*}\right)^{k-1}\left(1-\frac{1}{C^B}\left(\frac{w_{1}^{*}}{\tau_1}+\frac{w_{2}^{*}}{\tau_2}\right)^{B}\right),\notag\\
\mathcal{N}_{j} &=\frac{ B \left(w_{j}^{*}\right)^{2}}{\tau_{j}^{2}\left(C\right)^{B}}\left(\alpha \left(w_{j}^{*}\right)^{k-2}+\beta\right)\left(\frac{w_{1}^{*}}{\tau_{1}}+\frac{w_{2}^{*}}{\tau_{2}}\right)^{B-1},\notag\\
\mathcal{P}_{j}&=\frac{ B \left(w_{j}^{*}\right)^{2}}{\tau_{1}\tau_{2}\left(C\right)^{B}}\left(\alpha \left(w_{j}^{*}\right)^{k-2}+\beta\right)\left(\frac{w_{1}^{*}}{\tau_{1}}+\frac{w_{2}^{*}}{\tau_{2}}\right)^{B-1},
\end{align}
for $j=1,2$. Looking for exponential solutions, we obtain the characteristic equation for the linearised system \eqref{eq:linearmodelb_1reduced2} as
\begin{align}
\label{eq:linearmodel1a2kappachar}
\lambda^2 + \kappa a\lambda +\kappa b\lambda e^{-\lambda \tau_1} + \kappa^2 c e^{-\lambda \tau_1} + \kappa^2 d=0.
\end{align}
where, 
\begin{align}
\label{eq: abcd}
&a = \mathcal{M}_{1}+\mathcal{M}_{2}+\mathcal{N}_{2} \hspace{1ex}, \hspace{1ex} b=\mathcal{N}_{2},\notag\\
&c=\mathcal{M}_{1}\left(\mathcal{M}_{2}+\mathcal{N}_{2}\right), \hspace{2ex}d=\mathcal{N}_{1}\left(\mathcal{M}_{2}+\mathcal{N}_{2}\right)-
\mathcal{P}_{1}\mathcal{P}_{2}.
\end{align}
For system \eqref{eq:modelb_1reduced2} to be locally stable about the equilibrium $(w_{1}^{*},w_{2}^{*})$, all roots of the characteristic equation \eqref{eq:linearmodel1a2kappachar} should lie in the left half of the complex plane. It can be shown that, for negligibly small values of the non-dimensional parameter $\kappa$, the system is stable, \emph{i.e.}, all the roots would have negative real parts. However, as $\kappa$ is increased beyond a critical value, one pair of complex conjugate roots may cross over the imaginary axis, and hence have positive real parts. At this critical value the system would transit into an unstable region and have a pair of purely imaginary roots. To deduce this point, we substitute $\lambda=j\omega$ in \eqref{eq:linearmodel1a2kappachar} and separate real and imaginary parts to get
\begin{align*}
\omega^2 = \frac{\kappa^2(2c-a^2+b^2)}{2}\pm \frac{\kappa^2 \sqrt{(2c-a^2+b^2)^2-4(c^2-d^2)}}{2}.
\end{align*} 
\emph{Condition 1}: There exists only one positive value of $\omega^2$ if the following conditions hold
\begin{itemize}
\item[(i)]$(2c-a^2+b^2)>0,\text{and}\hspace{1mm} (2c-a^2+b^2)^2=4(c^2-d^2)$
\item[(ii)]$(2c-a^2+b^2)>0,\text{and}\hspace{1mm} c^2-d^2<0$
\end{itemize}
\emph{Condition 2}: There exists two positive value of $\omega^2$ if the following condition holds
\begin{align*}
(2c-a^2+b^2)>0,\text{and}\hspace{1mm} (c^2-d^2)>0.
\end{align*}
When Condition 1 is satisfied, the system transits from the locally stable regime to instability as $\kappa$ increases beyond a critical value, and never regains stability as $\kappa$ is further increased. On the contrary, when Condition 2 is satisfied, the system may undergo stability switches as $\kappa$ is increased~\cite{Cooke}. In the context of congestion control algorithms, the stability switch phenomenon is an undesirable dynamical feature. Further, we have observed in numerous packet-level simulations that Compound TCP does not exhibit stability switches. Hence, we focus only on the case when Condition 1 is satisfied, and only one positive root of $\omega^2$ exists. This implies that there exists a cross over frequency at which one pair of complex conjugate roots crosses over the imaginary axis, and is given by $\omega=\kappa A$, where
\begin{align*}
A=\sqrt{\frac{(2c-a^2+b^2)}{2}+\frac{\sqrt{(2c-a^2+b^2)^2-4(c^2-d^2)}}{2}}.
\end{align*}
The critical value of $\kappa$ denoted by $\kappa_{c}$, at which this transition occurs, is given by
\begin{align}
\kappa_{c}=\frac{1}{A\tau_1}\cos^{-1}\Big(\frac{A^2(d-ab)-cd}{b^2A^2+c^2}\Big).
\label{eq:kappa_c2}
\end{align}
\subsection*{Case II}
The schematic diagram of the topology is illustrated in Fig.~\ref{Figure.1}(b).
\subsubsection*{Scenario 1}With this assumption, the second model reduces to a single set of TCP flows, regulated by an edge router, and feeding into a core router. Observe that, the loss probabilities at both routers are the same. Hence, with the non-dimensional bifurcation parameter $\kappa$, system \eqref{eq:modelb_1} reduces to the following non-linear, first-order, time-delayed differential equation
\begin{align}
\frac{dw(t)}{dt} = \kappa\frac{w(t-\tau)}{\tau}\bigg(i\left(w(t)\right)\Big(1-p(w(t-\tau))\Big) - d\left((w(t)\right)q(w(t-\tau))\bigg),
\label{eq:modelb2_reduced1}
\end{align}
where $w(t)$ is the average window size of the TCP flows. Using the functional forms of loss probabilities given by \eqref{eq:case2_loss}, we obtain
\begin{align*}
p(w^{*})=2q(w^{*})=2\left(\frac{w^{*}}{C\tau}\right)^B.
\end{align*}
The critical value of $\kappa$, at which system \eqref{eq:modelb2_reduced1} loses its stability, satisfies the following equation 
\begin{align}
\kappa_{c}\alpha \left(w^{*}\right)^{k-1}\sqrt{B^{2}-\left(\left(k-2\right)\left(1-p(w^{*})\right)\right)^{2}} = \cos^{-1}\left(\frac{\left(k-2\right)\left(1-p(w^*)\right)}{B}\right).
\label{eq:condition_case2}
\end{align}
\subsubsection*{Scenario 2}With the introduction of the non-dimensional parameter $\kappa$, system \eqref{eq:modelb_2} reduces to
\begin{align}
\frac{dw_{1}(t)}{dt} =& \, \kappa \frac{w_{1}(t-\tau_{1})}{\tau_{1}}\bigg(i\left(w_{1}(t)\right)\Big(1-q_1(t,\tau_1,\tau_2)-q_2(t,\tau_1,\tau_2))\Big) - d\left((w_{1}(t)\right)\left(q_1(t,\tau_1,\tau_2)+q_2(t,\tau_1,\tau_2)\right)\bigg),\notag\\ 
\frac{dw_{2}(t)}{dt} =& \, \kappa \frac{w_{2}(t)}{\tau_{2}}\bigg(i\left(w_{2}(t)\right)\Big(1-q_1(t,\tau_1,\tau_2)-q_2(t,\tau_1,\tau_2))\Big) - d\left((w_{2}(t)\right)\left(q_1(t,\tau_1,\tau_2)+q_2(t,\tau_1,\tau_2)\right)\bigg).
\label{eq:modelb2_reduced2}
\end{align}
Linearising \eqref{eq:modelb2_reduced2} about its non-trivial equilibrium $(w_1^{*},w_2^{*})$, we obtain
\begin{align}
\label{eq:linearmodelb_2reduced2}
&\dot{u}_1(t) = -\kappa\left(\mathcal{M}_{1}u_{1}(t)+\mathcal{N}_{1}u_{1}(t-\tau_{1})+\mathcal{P}_{1}u_{2}(t)\right),\notag\\
&\dot{u}_2(t) = -\kappa\left(\big(\mathcal{M}_{2}+\mathcal{N}_{2}\big)u_{2}(t)+\mathcal{P}_{2}u_{1}(t-\tau_{1})\right),
\end{align}
where, for Compound TCP, the increase and decrease functions \eqref{eq:increasedecrease}, and the functional forms of the loss probabilities given by~\eqref{eq:case2_loss} yield the following coefficients 
\begin{align}
\label{eq:coefficients_2}
\mathcal{M}_{j} & = -\frac{\alpha}{\tau_{j}}\left(k-2\right) \left(w_{j}^{*}\right)^{k-1}\Bigg(1-\bigg(\frac{1}{C_{1}^{B_{1}}}\bigg(\frac{w_1^*}{\tau_1}+\frac{w_2^*}{\tau_2}\bigg)^{B_1} -\frac{1}{C_{2}^{B_{2}}}\bigg(\frac{w_1^*}{\tau_1}+\frac{w_2^*}{\tau_2}\bigg)^{B_2}\Bigg),\notag\\
\mathcal{N}_{j} & = \left(\alpha\left(w_{j}^{*}\right)^{k-1}+\beta w_{j}^{*}\right)\frac{w_{j}^{*}}{\tau_j}\Bigg(\frac{B_{1}}{C_{1}^{B_{1}}}\left(\frac{w_{1}^{*}}{\tau_{1}}+\frac{w_{2}^{*}}{\tau_{2}}\right)^{B_{1}-1} + \frac{B_2}{C_{2}^{B_{2}}}\left(\frac{w_{1}^{*}}{\tau_{1}}+\frac{w_{2}^{*}}{\tau_{2}}\right)^{B_{2}-1}\Bigg),\notag\\
\mathcal{P}_{j} & = \left(\alpha\left(w_{j}^{*}\right)^{k-1}+\beta w_{j}^{*}\right)\frac{w_{j}^{*}}{\tau_1\tau_2}\Bigg(\frac{B_{1}}{C_{1}^{B_{1}}}\left(\frac{w_{1}^{*}}{\tau_{1}}+\frac{w_{2}^{*}}{\tau_{2}}\right)^{B_{1}-1} +\frac{B_2}{C_{2}^{B_{2}}}\left(\frac{w_{1}^{*}}{\tau_{1}}+\frac{w_{2}^{*}}{\tau_{2}}\right)^{B_{2}-1}\Bigg),  
 \end{align}
for $j=1,2.$ Observe that, the linearised system \eqref{eq:linearmodelb_2reduced2} has a similar form as \eqref{eq:linearmodelb_1reduced2}. Hence, conducting a similar kind of analysis as done for system \eqref{eq:modelb2_reduced2}, we obtain the critical value of the non-dimensional parameter $\kappa$, as given by \eqref{eq:kappa_c2}.
\subsection*{Case III}
The schematic diagram for this topology is illustrated in Fig.~\ref{Figure.1}(c). 
\subsubsection*{Scenario 1}For Compound TCP in the small buffer regime, the critical value of $\kappa$, denoted by $\kappa_c$, at which system~\eqref{eq:modelb_3} transits into a locally unstable regime, satisfies the following condition:
\begin{align*}
\kappa_c\alpha \left(w^{*}\right)^{k-1}\sqrt{B^{2}
-(k-2)^2\left(1-\left(1+2^B\right)p(w^{*})\right)^{2}} < \cos^{-1}
\left(\frac{(k-2)\left(1-\left(1+2^B\right)p(w^*)\right)}{B}\right).
\end{align*}
\begin{figure*}
\begin{center}
  \begin{minipage}[b]{0.45\textwidth}
  \label{fig:alpha_kappa}
  \psfrag{b}{Non-dimensional parameter, $\kappa$}
  \psfrag{T}{Protocol parameter, $\alpha$}
  \psfrag{h}{\begin{scriptsize}Hopf condition\end{scriptsize}}
  \psfrag{0}{\begin{scriptsize}$0$\end{scriptsize}}
  \psfrag{2}{\begin{scriptsize}$2$\end{scriptsize}}
  \psfrag{1}{\begin{scriptsize}$1$\end{scriptsize}}
\psfrag{0.3}{\begin{scriptsize}$0.3$\end{scriptsize}}  
     \includegraphics[width=2in,height=3in,angle=270]{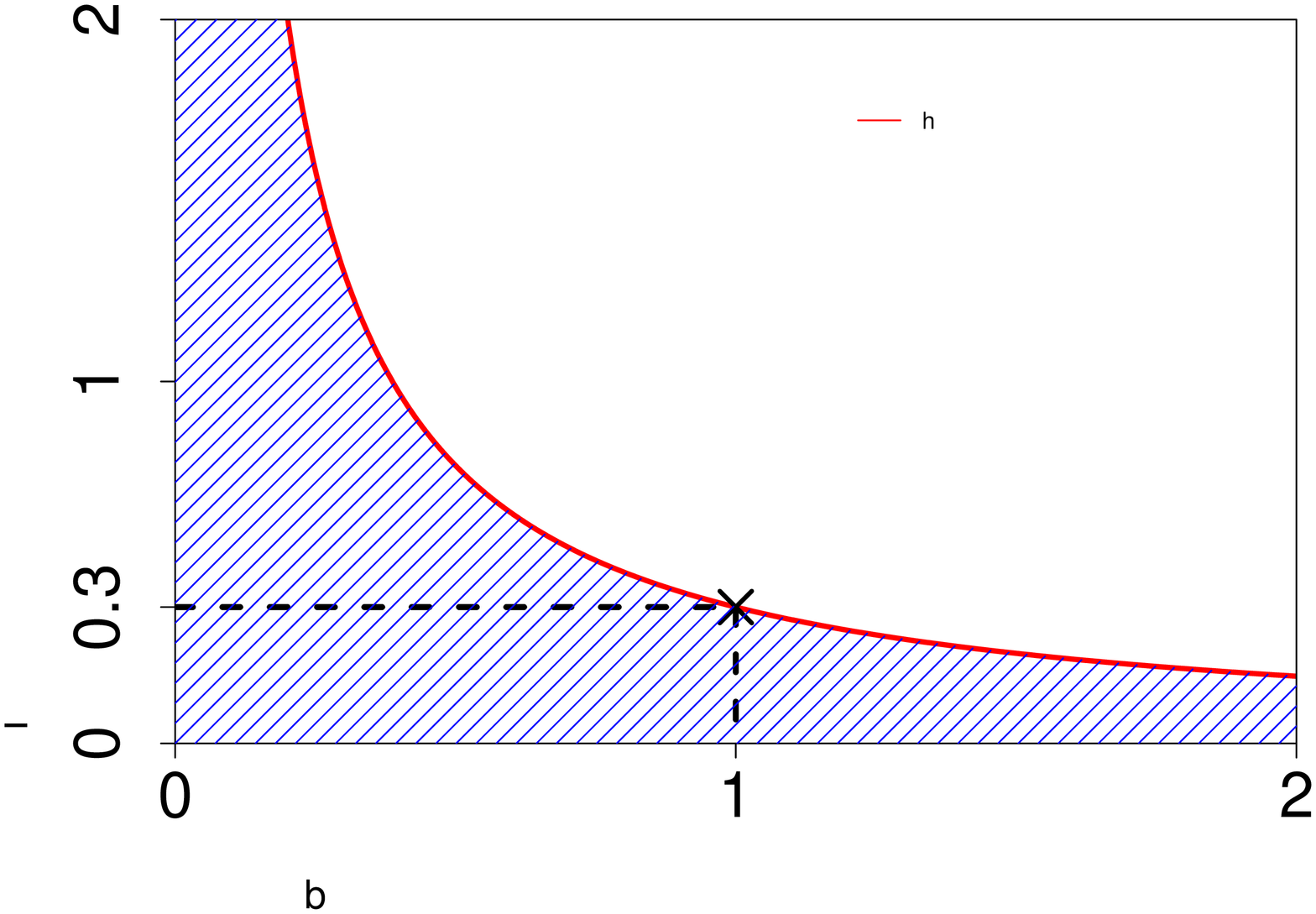}
  
  \end{minipage}
 \hspace*{5mm}
  \begin{minipage}[b]{0.45\textwidth}
  \label{fig:buffer_kappa}
  \psfrag{b}{Non-dimensional parameter, $\kappa$}
  \psfrag{T}{Core router buffer, $B$}
  \psfrag{h}{\begin{scriptsize}Hopf condition\end{scriptsize}}
  \psfrag{0.9}{\begin{scriptsize}$0.9$\end{scriptsize}}
  \psfrag{1.2}{\begin{scriptsize}$1.2$\end{scriptsize}}
  \psfrag{1.5}{\begin{scriptsize}$1.5$\end{scriptsize}}
   \psfrag{16}{\begin{scriptsize}$16$\end{scriptsize}}
  \psfrag{22}{\begin{scriptsize}$22$\end{scriptsize}}
  \psfrag{28}{\begin{scriptsize}$28$\end{scriptsize}}
    \includegraphics[width=2in,height=3in,angle=270]{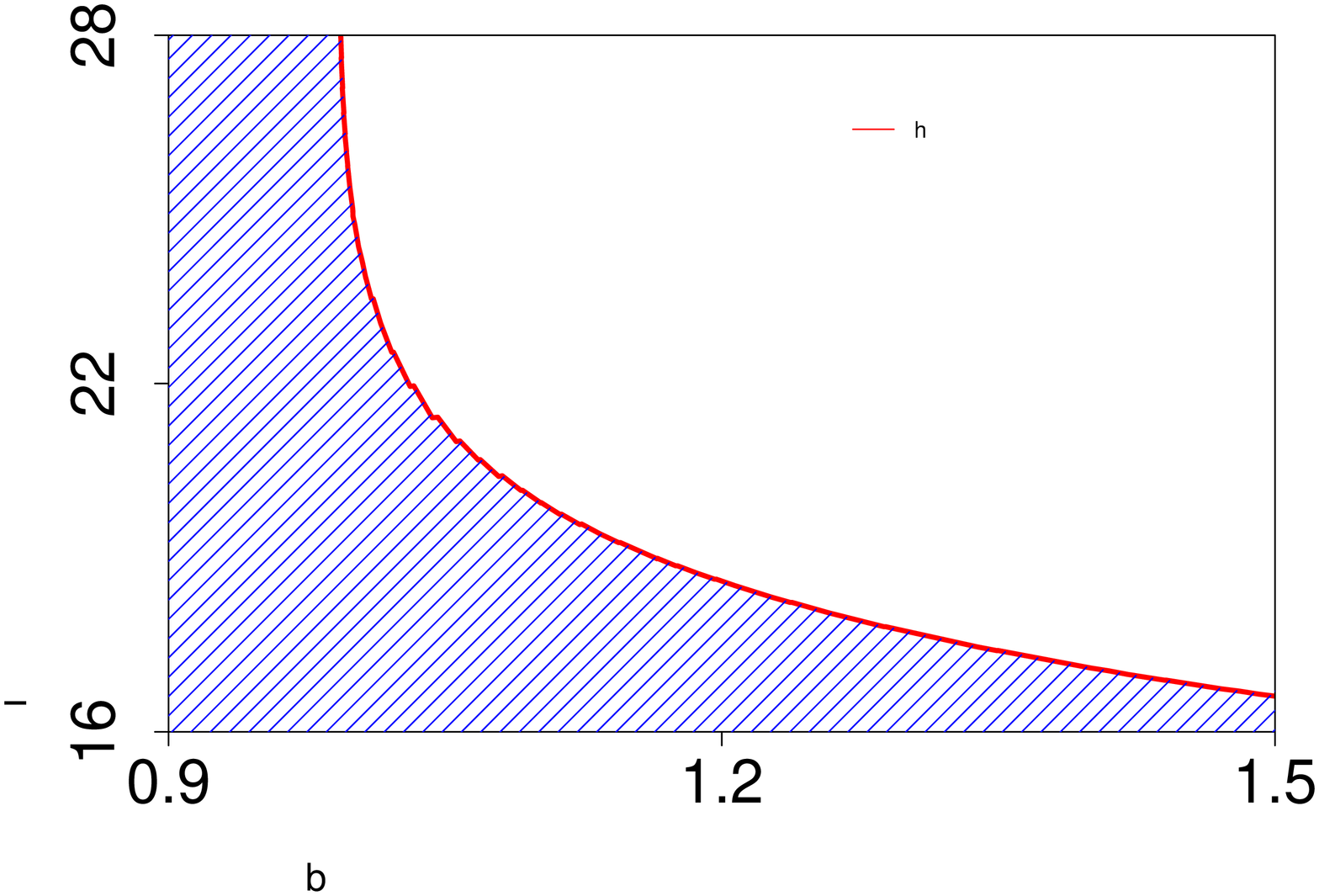}
    \end{minipage}
  
     \caption {\emph{Stability chart.} Hopf condition for~\eqref{eq:modelb_3} with Compound TCP in the small buffer regime with respect to two sets of parameters: (a) the non-dimensional parameter $\kappa$, and the protocol parameter $\alpha$, (b) the non-dimensional parameter $\kappa$, and the buffer size at the core router $B$. The shaded region below the Hopf condition curve represents the stable region.}
  \label{fig:charts}
\end{center}
\end{figure*}
\subsubsection*{Scenario 2}With the introduction of the non-dimensional parameter $\kappa$, \eqref{eq:modelb_3} reduces to 
\begin{align}
\frac{dw_{1}(t)}{dt} & = \kappa\frac{w_{1}(t-\tau_{1})}{\tau_{1}}\bigg(i\left(w_{1}(t)\right)\Big(1-p_1(t-\tau_1)-q(t,\tau_1,\tau_2))\Big) - d\left((w_{1}(t)\right)\left(p_1(t-\tau_1)+q(t,\tau_1,\tau_2)\right)\bigg),\notag\\ 
\frac{dw_{2}(t)}{dt} & = \kappa \frac{w_{2}(t)}{\tau_{2}}\bigg(i\left(w_{j}(t)\right)\Big(1-p_2(t)-q(t,\tau_1,\tau_2))\Big) - d\left((w_{j}(t)\right)\left(p_2(t)+q(t,\tau_1,\tau_2)\right)\bigg).
\label{eq:modelb3_reduced2}
\end{align}
Linearising \eqref{eq:modelb3_reduced2} about its equilibrium $(w_{1}^{*},w_{2}^{*})$, we obtain
\begin{align}
\label{eq:linearmodelb_3reduced2}
&\dot{u}_1(t) = -\kappa\left(\mathcal{M}_{1}u_{1}(t)+\mathcal{N}_{1}u_{1}(t-\tau_{1})+\mathcal{P}_{1}u_{2}(t)\right),\notag\\
&\dot{u}_2(t) = -\kappa\left(\big(\mathcal{M}_{2}+\mathcal{N}_{2}\big)u_{2}(t)+\mathcal{P}_{2}u_{1}(t-\tau_{1})\right),
\end{align}
where, for Compound TCP, and the functional forms of the loss probabilities given by~\eqref{eq:case3_loss} yield the following coefficients
\begin{align}
\label{eq:coefficients3}
\mathcal{M}_{j} & = -\frac{\alpha}{\tau_{j}}\left(k-2\right) \left(w_{j}^{*}\right)^{k-1}\Bigg(1-\bigg(\frac{w_j^*}{C_j\tau_j}\bigg)^{B_j} -\frac{1}{C^B}\bigg(\frac{w_1^*}{\tau_1}+\frac{w_2^*}{\tau_2}\bigg)^{B}\Bigg),\notag\\
\mathcal{N}_{j} & = \left(\alpha\left(w_{j}^{*}\right)^{k-1}+\beta w_{j}^{*}\right)\Bigg(\frac{B_{j}}{\tau_j}\left(\frac{w_{j}^{*}}{C_j \tau_j}\right)^{B_j} +\frac{B \left(w_{j}^{*}\right)^{2}}{C^{B}\tau_{j}^{2}}\left(\frac{w_{1}^{*}}{\tau_{1}}+\frac{w_{2}^{*}}{\tau_{2}}\right)^{B-1}\Bigg),\notag\\
\mathcal{P}_{j} & = \left(\alpha\left(w_{j}^{*}\right)^{k-1}+\beta w_{j}^{*}\right)\frac{B w_{j}^{*}}{\tau_{1}\tau_{2}\left(C\right)^{B}}\left(\frac{w_{1}^{*}}{\tau_{1}}+\frac{w_{2}^{*}}{\tau_{2}}\right)^{B-1}, 
\end{align}
for $j=1,2.$ Note that, the linearised system \eqref{eq:linearmodelb_3reduced2} has a similar form as \eqref{eq:linearmodelb_1reduced2}. Hence, a similar kind of local stability analysis would yield the condition on the critical value of the non-dimensional parameter $\kappa$, and the protocol parameters as given by \eqref{eq:kappa_c2}.
\\
\indent For all three scenarios, with the simplifying assumptions, the conditions derived above essentially capture the interdependence among the non-dimensional parameter $\kappa$, and the system parameters to ensure local stability. Observe that, the loss of local stability can be studied with respect to any system parameter. However, we prefer
to choose an exogenous parameter as the bifurcation parameter, to aid our analysis. It can be explicitly shown that, for all the above cases, the system loses local stability via a Hopf bifurcation~\cite{Hassard} if the conditions derived above get violated. We prove this by verifying that the transversality condition of the Hopf spectrum \cite{Hassard}. To verify this, we show that, $\mathrm{Re}(\mathrm{d}\lambda/\mathrm{d}\kappa)\neq0$ at $\kappa=\kappa_c$. In particular, we prove that, $\mathrm{Re}(\mathrm{d}\lambda/\mathrm{d}\kappa)>0$ at $\kappa=\kappa_c$. This implies that, one pair of complex conjugate roots crosses over the imaginary axis from the left half of the complex plane to the right half. Thus, the system undergoes a Hopf bifurcation at $\kappa=\kappa_c$. Hence, $\kappa<\kappa_c$ is a necessary and sufficient condition for local stability, for all the three scenarios.\\     
\indent Observe that, deriving a necessary and sufficient condition with heterogeneous network parameters, and different round trip times is analytically complex, for all three scenarios discussed earlier. Hence, we numerically illustrate through DDE-BIFTOOL version 2.03 ~\cite{DDE1},~\cite{DDE2}, that system \eqref{eq:modelb_3} undergoes a Hopf bifurcation if the non-dimensional parameter $\kappa$ is varied beyond a certain critical value. We fix the protocol parameters as follows: $\alpha=0.3$, $\beta=0.5$ and $k=0.75$. Since, we mainly focus on small buffer regime, the buffer sizes of the routers are fixed as: $B_1=10, B_2=15$, and $B=25$. We fix the remaining network parameters as: $C_1=C_2=100$, $C=180$, $\tau_1=1$ and $\tau_2=2$. Now, we vary the non-dimensional parameter $\kappa$ in the range $[0,2]$ and observe that the system undergoes a Hopf bifurcation at $\kappa_c=1$. At this point, the system has one pair of complex conjugate roots on the imaginary axis. Consequently, the system dynamics exhibit limit cycles at $\kappa_c=1$.
 
\subsubsection*{Stability charts}
\begin{figure}
\begin{center}
 \psfrag{T}{$k$}
  \psfrag{b}{$\alpha$}
  \psfrag{0.1}{\begin{scriptsize}$0.1$\end{scriptsize}}
  \psfrag{1.0}{\begin{scriptsize}$1.0$\end{scriptsize}}
  \psfrag{0.5}{\begin{scriptsize}$0.5$\end{scriptsize}}
      \psfrag{SSSSSSSSSSSSSS}{\hspace{2mm}\begin{scriptsize}Stable Region\end{scriptsize}}
  \psfrag{h}{\begin{scriptsize}Hopf Condition\end{scriptsize}}
  \includegraphics[width=1.8in,height=2.6in,angle=270]{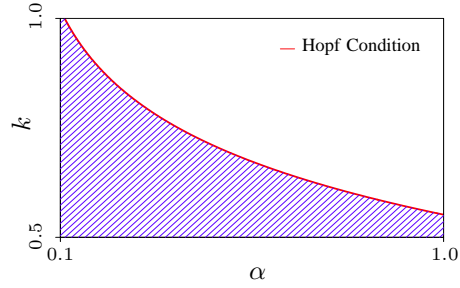}
  \caption{\emph{Stability chart.} Hopf condition for system \eqref{eq:modelb_3}. The shaded region below the curve denotes the stable region. It is evident that there exists a trade-off between the increase protocol parameters $\alpha$ and $k$. As $\alpha$ increases, $k$ has to be decreased to ensure system stability. }\label{fig:stability_chart}\vspace{-8mm}
  \end{center}
\end{figure}

To obtain insights about the system behaviour at the stability boundary, we now demonstrate some stability charts for system \eqref{eq:modelb_3}. Fig.~\ref{fig:charts}~(a) represents the Hopf condition for system \eqref{eq:modelb_3} in the two parameter space: the non-dimensional parameter $\kappa$, and the protocol parameter $\alpha$. Observe that, if $\kappa$ is increased, $\alpha$ would have to reduce to ensure stability. Fig.~\ref{fig:charts}~(b) illustrates the Hopf condition in the two parameter space: the non-dimensional parameter $\kappa$, and the buffer size at the core router $B$. Observe that, if $\kappa$ is increased, keeping other system parameters fixed, $B$ would have to be decreased accordingly to ensure stability of system~\eqref{eq:modelb_3}. Fig.~\ref{fig:stability_chart} characterises the stability boundary of system~\eqref{eq:modelb_3} with respect to the increase protocol parameters $\alpha$ and $k$. It is evident that, there exists a trade-off between the increase parameters to ensure stability. Hence, we conclude that, both protocol parameters, and network parameters, need to be co-designed carefully to maintain stability of system~\eqref{eq:modelb_3}. If these Hopf conditions get violated, the system would lose stability leading to the emergency of limit cycles in the system dynamics. In the next section, we provide a detailed analytical framework to characterise the \emph{type} of Hopf bifurcation and the asymptotic \emph{orbital stability} of the emergent limit cycles, for system \eqref{eq:modelb_3}.
\section{Hopf Bifurcation Analysis}
\label{sec:hopf}
We have seen that, variation in the exogenous parameter $\kappa $ induces instability in system. Instability in the system could be induced by any of the system parameters. This loss of stability occurs via a Hopf bifurcation which results in limit cycles in the system dynamics which in turn leads to deterministic oscillations in the queue size. Consequently, this results in the overall degradation of the system performance because of loss in link utilisation. To that end, it becomes imperative to study the type of bifurcation and the stability of these emergent limit cycles to a greater detail. 

Note that, we have motivated the exogenous, non-dimensional parameter $\kappa$ as the bifurcation parameter. This enables us to capture the effect of the different system parameters on the system stability in a unified manner. The Hopf bifurcation analysis enables us to analyse the system dynamics in its locally unstable regime, in the neighbourhood of the Hopf condition. Using Poincar\'{e} normal forms and the center manifold theory, we present an analytical framework to determine the \emph{type} of the Hopf bifurcation and the orbital stability of the emergent limit cycles. Our analysis closely follows the analysis presented in \cite{Hassard,Kuznetsov,Gaurav}.

Let $\kappa=\kappa_c+\mu$, where $\mu \in \mathbb{R}.$ Observe that, the system undergoes a Hopf bifurcation at $\mu=0$, where $\kappa=\kappa_c$. We can now consider $\mu$ as the bifurcation parameter. An incremental change in $\kappa$ from $\kappa_c$ to $\kappa_c+\mu$ where $\mu>0$, pushes the system to its locally unstable regime.

\emph{Step 1}: Using Taylor series expansion, we segregate the right hand side of \eqref{eq:modelb_3} into linear and non-linear parts. We then cast this equation into the standard form of an operator differential equation.

\emph{Step 2}: At the critical value of the bifurcation parameter, \emph{i.e.} at $\mu=0$, the system has exactly one pair of purely imaginary eigenvalues with non-zero angular velocity. The linear eigenspace spanned by the  eigenvectors corresponding these eigenvalues is called the critical eigenspace. The center manifold theorem~\cite[Chapter $5$, Theorem $5.1.$]{Kuznetsov} guarantees the existence of a locally invariant $2-$dimensional manifold which is tangent to the critical eigenspace at the equilibrium of the system.

\emph{Step 3}: Next, we project the system onto its critical eigenspace and its complement at the critical value of the bifurcation parameter. This enables us to capture the dynamics of the system on the center manifold, with the help of an ordinary differential equation in a single complex variable.

\emph{Step 4}: Finally, using Poincar\'{e} normal forms, we evaluate the lyapunov coefficient and the floquet exponent, which characterise the type of the Hopf bifurcation and the asymptotic orbital stability of the emergent limit cycles respectively.

Suppose $(w_1^*,w_2^*)$ is an equilibrium for \eqref{eq:modelb_3}. Let $u_1(t)=w_1(t)-w_1^{\ast}$ and $u_2(t)=w_2(t)-w_2^{\ast}$ be small perturbations about the equilibrium. Thus, a Taylor series expansion of \eqref{eq:modelb_3} about its equilibrium $\left(w_1^{\ast},w_2^{\ast}\right)$ is as follows
\begin{align}
\dot{u_1}(t) = & \, \kappa \Big( \xi_{a} u_1(t)+ \xi_{b} u_1(t-\tau_ 1)+\xi_{d}u_2(t-\tau_2) \Big) +\kappa\Big(\xi_{aa}u_1^2(t)+\xi_{bb}u_1^2(t-\tau_1)  +\xi_{dd}u_2^2(t-\tau_2) \notag \\
& + \xi_{ab}u_1(t)u_1(t-\tau_1)+\xi_{ad}u_1(t)u_2(t-\tau_2) + \xi_{bd}u_1(t-\tau_1)u_2(t-\tau_2)\Big) +\kappa\Big(\xi_{aaa}u_1^3(t) \notag \\
& +\xi_{bbb}u_1^3(t-\tau_1)+\xi_{ddd}u_2^3(t-\tau_2) + \xi_{aab}u_1^2(t)u_1(t-\tau_1) +\xi_{aad}u_1^2(t)u_2(t-\tau_2) + \xi_{abb}u_1(t)u_1^2(t-\tau_1) \notag \\
& +\xi_{bbd}u_1^2(t-\tau_1)u_2(t-\tau_2)+
\xi_{add}u_1(t)u_2^2(t-\tau_2)+\xi_{bdd}u_1(t-\tau_1)u_2^2(t-\tau_2) \notag \\
& + \xi_{acd}u_1(t)u_1(t-\tau_1)u_2(t-\tau_2)\Big),\notag \\
\dot{u_2}(t) = & \, \kappa \Big( \chi_{c} u_2(t)+ \chi_{d} u_2(t-\tau_ 2)+\chi_{b}u_1(t-\tau_1) \Big) +\kappa \Big(\chi_{cc}u_2^2(t)+\chi_{dd}u_2^2(t-\tau_2) + \chi_{bb}u_1^2(t-\tau_1)\notag\\
& + \chi_{cd}u_2(t)u_2(t-\tau_2)+\chi_{bc}u_1(t-\tau_1)u_2(t) + \chi_{bd}u_1(t-\tau_1)u_2(t-\tau_2)\Big) +\kappa\Big(\chi_{ccc}u_2^3(t) \notag \\
& +\chi_{ddd}u_2^3(t-\tau_2)+\chi_{bbb}u_1^3(t-\tau_1) + \chi_{ccd}u_2^2(t)u_2(t-\tau_2)+\chi_{bcc}u_1(t-\tau_1)u_2^2(t) + \chi_{cdd}u_2(t)u_2^2(t-\tau_2) \notag \\
& +\chi_{bdd}u_1(t-\tau_1)u_2^2(t-\tau_2) +\chi_{bbc}u_1^2(t-\tau_1)u_2(t)+\chi_{bbd}u_1^2(t-\tau_1)u_2(t-\tau_2)\notag\\
& + \chi_{bcd}u_1(t-\tau_1)u_2(t)u_2(t-\tau_2)\Big).\label{eq:TaylorSeries}
\end{align}
The Taylor series coefficients are given in Table \ref{table:coefficients}. Using the notation $\mathbf{u} = [u_1 \hspace{4pt} u_2]^T$, we reduce equation \eqref{eq:TaylorSeries} to the following form
\begin{align}
\dot{\mathbf{u}}(t) = \mathcal{L}_{\mu}\mathbf{u}_{t} + \mathcal{F}(\mathbf{u}_{t},\mu),\label{eq:nonlinear}
\end{align}
where $t>0, \mu \in \mathbb{R}$. For $\tau>0$,we define  
\begin{align*}
\mathbf{u}_{t}(\theta)=\mathbf{u}(t+\theta), \hspace{2ex} \mathbf{u}_{t}:[-\tau,0] \rightarrow \mathbb{R}^2, \hspace{2ex} \theta \in [-\tau,0].
\end{align*}
For this model, without loss of generality, we assume that $\tau_1>\tau_2$. $\mathcal{L}: C[-\tau_1,0] \rightarrow \mathbb{R}^{2}$ denotes a family of continuous and bounded functions parametrised by $\mu$. Here, $C[a,b]$ denotes the set of all continuous functions on the interval $[a,b]$. The operator $\mathcal{F}:C[-\tau_1,0] \rightarrow \mathbb{R}^{2}$ consists of the non-linear terms. Further, we assume that $\mathcal{F}$ is analytic and both $\mathcal{L}$ and $\mathcal{F}$ depend analytically on the bifurcation parameter $\mu$ for small $|\mu|$. The linear operator is
\begin{align}
\label{eq:linearop} 
\mathcal{L}_{\mu} \mathbf{u}_{t} = \kappa \begin{bmatrix}
                                   \xi_{a}u_1(t)+\xi_{b}u_1(t-\tau_1)    & \xi_{d}u_2(t-\tau_2)\\
                                   \chi_{b}u_1(t-\tau_1)  & \chi_{c}u_2(t)+\chi_{d}u_2(t-\tau_2) \\
                                 \end{bmatrix}.          
\end{align}
We now cast equation \eqref{eq:nonlinear} into the following standard form of an operator differential equation,
\begin{align}
\label{eq:system}
\dot{\mathbf{u}} = \mathcal{A}(\mu)\mathbf{u}_{t}+\mathcal{R}\mathbf{u}_{t}.
\end{align} 
Note that, \eqref{eq:system} has $\mathbf{u}_{t}$ rather than both $\mathbf{u}_{t}$ and $\mathbf{u}$. Now, using the Riesz representation theorem~\cite[Chapter $6,$ Theorem $6.19.$]{Rudin}, we transform the linear problem $(\mathrm{d}/\mathrm{d}t)\mathbf{u}(t)=\mathcal{L}_{\mu}\mathbf{u}_{t}$. The Riesz representation theorem guarantees the existence of an $2\times 2$ matrix-valued measure $\boldsymbol{\eta}(\cdot,\mu):[-\tau_1,0]\rightarrow\mathbb{R}^{4}$, such that each component of $\boldsymbol{\eta}$ has bounded variation and for all  $\phi \in C[-\tau_1,0],$
\begin{align*}
 \mathcal{L}_{\mu} \phi =& \int_{\theta = -\tau_1}^{0} \mathrm{d} \boldsymbol{\eta}(\theta, \mu) \phi(\theta).
 \end{align*}
 In particular, we have
 \begin{align*}
 \mathcal{L}_{\mu} \mathbf{u}_{t} =& \int_{\theta = -\tau_1}^{0} \mathrm{d} \boldsymbol{\eta}(\theta, \mu) \mathbf{u}_{t}(\theta).\notag
 \end{align*}
Observe that, for system \eqref{eq:modelb_3}, the matrix $\mathrm{d}\boldsymbol{\eta}$ is  \begin{align}
 \label{eq:etamatrix}
 \mathrm{d} \boldsymbol{\eta}(\theta,\mu) = \kappa \begin{bmatrix}
                                   \xi_{a}\delta(\theta)+\xi_{b}\delta(\theta+\tau_1)    & \xi_{d}\delta(\theta+\tau_2)\\
                                   \chi_{b}\delta(\theta+\tau_1)  & \chi_{c}\delta(\theta)+\chi_{d}\delta(\theta+\tau_2) \\
                                 \end{bmatrix}\mathrm{d}\theta.
 \end{align}
Here, $\delta(\cdot)$ is the Dirac-delta measure. Let  $C^{1}[-\tau_1,0]$ denotes the space of all functions defined on $[-\tau_1,0]$, with continuous first derivatives. For $\phi \in C^{1}[-\tau_1,0],$ we then define the following linear and non-linear operators 
\begin{align}
\mathcal{A}(\mu) \mathbf{u}_{t}(\theta) &= \left\lbrace \begin{array}{lr}
\frac{\mathrm{d} \mathbf{u}_{t}(\theta)}{\mathrm{d}\theta}, & \theta \in [-\tau_1,0), \\
\mathcal{L}_{\mu}\mathbf{u}_{t}, & \theta = 0.
 \end{array} \right. \notag\\
\mathcal{R}\mathbf{u}_{t}(\theta) &= \left \{ \begin{array}{cl}
0, & \theta \in [-\tau_1,0),\\
\mathcal{F}(\mathbf{u}_{t},\mu), & \theta = 0.
\end{array} \right.  \label{eq:DLEoperatorR}
\end{align}
Note that, $\mathrm{d}\mathbf{u}_{t}/\mathrm{d}\theta \equiv \mathrm{d}\mathbf{u}_{t}/\mathrm{d}t$. Hence, equation \eqref{eq:nonlinear} can be transformed into \eqref{eq:system}. Further, recall that, $\kappa=\kappa_{c}+\mu$, and the system undergoes bifurcation at the critical point $\mu=0$. Hence, we fix $\mu=0$ to perform the necessary analysis at the point of bifurcation. At $\mu=0$, the system has a pair of complex eigenvalues on the imaginary axis: $\lambda=\pm i\omega_0$, where $\omega_0>0.$ Let $\mathbf{q}(\theta)$ denote the eigenvector for $\mathcal{A}(0)$ corresponding to the eigenvalue $\lambda(0)=i\omega_0$. We assume that $\mathbf{q}(\theta)$ has a form as  
\begin{align*}
 \mathbf{q}(\theta) &= \begin{bmatrix}1 & \phi_1 \end{bmatrix}^{T}e^{i\omega_0\theta}.
\end{align*}
Now, using
\begin{equation*}
\mathcal{A}(0)\mathbf{q}(\theta) = i\omega_{0}\mathbf{q}(\theta),
\end{equation*}
we obtain $\phi_1$ as
\begin{align*}
\phi_1=\frac{-\kappa \chi_{b}e^{-i\omega_0 \tau_1}}{\kappa \chi_{c}+\kappa \chi_{d}e^{-i\omega_0 \tau_2}-i\omega_0}.
\end{align*}
We now define the following adjoint operator
\begin{equation*}
 \mathcal{A}^{\ast}(\mu)\alpha(s) = \begin{cases} \begin{array}{ll}-\frac{\mathrm{d}\alpha(s)}{\mathrm{d}s},&s \in (0,\tau_1],\\\int_{t = -\tau_1}^{0}\mathrm{d}\boldsymbol{\eta}^{T}(t,0)\alpha(-t),&s=0.\end{array}\end{cases}
\end{equation*}
where $\boldsymbol{\eta}^{T}$ denotes the transpose of $\boldsymbol{\eta}$. Observe that, the domains of $\mathcal{A}$ and $\mathcal{A}^{\ast}$ are $C^{1}[-\tau_1,0]$ and $C^{1}[0,\tau_1]$ respectively. Then, $\bar{\lambda}(0)=-i\omega_0$ is an eigenvalue of $\mathcal{A}^\ast$ and for some non-zero vector $\mathbf{p}$, we have
\begin{equation}
\label{eq:A*}
\mathcal{A}^{\ast}(0)\boldsymbol{\mathrm{p}}(\zeta) = -i\omega_{0}\boldsymbol{\mathrm{p}}(\zeta).
\end{equation}
We consider $\mathbf{p}(\theta)$ to have the following form: 
\begin{align*}
 \mathbf{p}(\theta) &= D\begin{bmatrix}\phi_2 & 1\end{bmatrix}^{T}e^{i\omega_0\theta}.
 \end{align*}
Using \eqref{eq:A*}, we obtain $\phi_2$ as
 \begin{align*}
 \phi_2=\frac{-\kappa \chi_{b}e^{i\omega_0\tau_1}}{\kappa \xi_{a}+\kappa \xi_{b}e^{i\omega_0\tau_1}+i\omega_0}.
 \end{align*}
 Let us define the inner-product of the functions $\boldsymbol{\psi} \in C[0,\tau_1]$ and $\boldsymbol{\phi} \in C[-\tau_1,0]$ as 
\begin{align}
\left\langle \boldsymbol{\psi}, \boldsymbol{\phi} \right\rangle =& \boldsymbol{\bar{\psi}}(0) \boldsymbol{\phi}(0) - \int_{\theta = -\tau_1}^{0} \int_{\zeta=0}^{\theta} \overline{\boldsymbol{\psi}}^{T}(\zeta - \theta) \mathrm{d}\boldsymbol{\eta}(\theta, \mu) \boldsymbol{\phi}(\zeta)\mathrm{d}\zeta. \label{eq:InnerProductDefn}
\end{align}
Using the above definition of inner product, we can easily verify that the eigenvectors $\mathbf{p}$ and $\mathbf{q}$ satisfy the conditions $\left\langle \mathbf{p},\mathbf{q}\right\rangle = 1$ and $\left\langle \mathbf{p},\mathbf{\bar{q}}\right\rangle = 0$ when
\begin{align}
D &=\Big(\phi_{2}\left(1+\kappa \xi_{b}\tau_1 e^{-i\omega_0\tau_1}+\kappa \xi_{d}\phi_1 \tau_2 e^{-i\omega_0 \tau_2}\right)+\phi_{1}\left(1+\kappa \chi_{d}\tau_2 e^{-i\omega_0\tau_2}+
\kappa \chi_{b}\tau_1 e^{-i\omega_0 \tau_1}\right)\Big)^{-1}.
\label{eq:exprn_B}
\end{align}
The critical eigenspace corresponding to the pair of eigenvalues $\pm i\omega_0$, denoted by $T_c$, is now $2-$dimensional and is spanned by $\lbrace\mathrm{Re}\hspace{0.25ex}\boldsymbol{\mathrm{q}}, \hspace{0.5ex}\mathrm{Im}\hspace{0.25ex}\boldsymbol{\mathrm{q}}\rbrace,$ where $\mathrm{Re}\hspace{0.25ex}\boldsymbol{\mathrm{q}}$ and $\mathrm{Im}\hspace{0.25ex}\boldsymbol{\mathrm{q}}$ denote the real and imaginary parts of $\boldsymbol{\mathrm{q}}$ respectively. Further, we denote the complement of the critical eigenspace $T_c$ as $T_{su}$. We now project system \eqref{eq:system} onto $T_c$ and $T_{su}$. For $\mathbf{u}_{t}$, a solution of \eqref{eq:system}  at $\mu = 0$, define
\begin{align}
  z(t) = \langle \mathbf{p},\mathbf{u}_{t}\rangle ,\hspace{3mm} \text{and} \hspace{3mm} \mathbf{w}(t,\theta) &= \mathbf{u}_{t}(\theta)- 2 \text{Re}\big(z(t)\mathbf{q}(\theta)\big).
 \label{eq:eigenbasis}
\end{align}

Recall that, the center manifold, $C_0$ is tangent to the critical eigenspace at the equilibrium. The representation of the center manifold is  
\begin{align}
\mathbf{w}(t,\theta) &= \mathbf{w}\big(z(t),\bar{z}(t),\theta\big),\hspace{1ex}\text{where}\notag\\
\mathbf{w}(z,\bar{z},\theta) &= \mathbf{w}_{20}(\theta)\frac{z^2}{2} + \mathbf{w}_{11}(\theta)z\bar{z} + \mathbf{w}_{02}(\theta)\frac{\bar{z}^2}{2} + \cdots.\label{eq:vector_w}
\end{align}
Here, $\mathbf{w}_{ij}(\theta),\hspace{0.5ex}\text{for all}\hspace{0.5ex}i,j\in \lbrace 0,1,2\rbrace$ is a two dimensional vector given as 
\begin{align*}
\mathbf{w}_{ij}(\theta)=\begin{bmatrix}
w_{ij1}(\theta) & w_{ij2}(\theta)
\end{bmatrix}^{T}.
\end{align*}
We observe that, $z$ and $\bar{z}$ are the local coordinates on the manifold $C_{0}$ in the direction of the eigenvectors $\mathbf{p}$ and $\bar{\mathbf{p}}$ respectively. Further, note that the existence of the center manifold $C_0$ ensures that equation \eqref{eq:system} can now be reduced to an ordinary differential equation for a single complex variable $z$ on $C_0$. At $\mu = 0$, in the coordinates $\eqref{eq:eigenbasis},$ the dynamics of $z$ can be represented as 
\begin{align}
   \dot{z}(t) &= \langle \mathbf{p},\mathcal{A}\mathbf{u_{t}} + \mathcal{R}\mathbf{u_{t}}\rangle \notag\\
   & = i\omega_{0}z(t) + \bar{\mathbf{p}}(0)\cdot \mathcal{F}\Big(\mathbf{w}(z,\bar{z},\theta)+ 2\text{Re}\big(z(t)\mathbf{q}(\theta)\big)\Big)\notag\\
   & = i\omega_{0}z(t) + \bar{\mathbf{p}}(0)\cdot \mathcal{F}_{0}(z,\bar{z}) \notag\\
   & = i\omega_{0}z(t) + g(z,\bar{z}).\label{eq:derivative_z}
\end{align}
Now, we can expand the function $g(z,\bar{z})$ in powers of $z$ and $\bar{z}$ as 
\begin{align}
g(z,\bar{z}) =&\,\, g_{20}\frac{z^{2}}{2} + g_{11}z \bar{z} + g_{02}\frac{\bar{z}^{2}}{2}+ g_{21}\frac{z^{2}\bar{z}}{2}+ \cdots. \label{eq:vector_g}
\end{align}
We now need to determine the coefficients $\mathbf{w}_{11}(\theta)$, $\mathbf{w}_{20}(\theta)$, $\mathbf{w}_{02}(\theta)$ in equation \eqref{eq:vector_w} to solve the differential equation \eqref{eq:derivative_z} for $z$. Following \cite{Hassard} we can write $\dot{\mathbf{w}} = \dot{\mathbf{u}}_{t} - \dot{z}\mathbf{q} - \dot{\bar{z}}\bar{\mathbf{q}},$ and using \eqref{eq:system} and \eqref{eq:derivative_z} we obtain
\begin{equation*}
\dot{\mathbf{w}} = \begin{cases}
               \begin{array}{l l}
                \mathcal{A}\mathbf{w} - 2\text{Re}\big(\bar{\mathbf{p}}(0)\cdot\mathcal{F}_0\mathbf{q}(\theta)\big),& \theta \in [-\tau_1,0),\\
                \mathcal{A}\mathbf{w} - 2\text{Re}\big(\bar{\mathbf{p}}(0)\cdot\mathcal{F}_0\mathbf{q}(0)\big)+\mathcal{F}_0, & \theta = 0,
               \end{array}
              \end{cases}
\end{equation*}
which, using \eqref{eq:vector_w}, can be rewritten as 
\begin{equation}
 \dot{w} = \mathcal{A}\mathbf{w} + \mathbf{H}(z,\bar{z},\theta).\label{eq:derivative_w}
\end{equation}
Here, the function $\mathbf{H}(z,\bar{z},\theta)$ can be expanded in powers of $z$ and $\bar{z}$ as 
\begin{equation}
\mathbf{H}(z,\bar{z},\theta) = \mathbf{H}_{20}(\theta)\frac{z^{2}}{2} + \mathbf{H}_{11}(\theta)z\bar{z} + \mathbf{H}_{02}(\theta)\frac{\bar{z}^{2}}{2} + \cdots.\label{eq:vector_H}
\end{equation}
Here, $\mathbf{H}_{ij}(\theta),\hspace{0.5ex}\text{for all}\hspace{0.5ex}i,j\in \lbrace 0,1,2\rbrace$ is a two dimensional vector given as 
\begin{align*}
\mathbf{H}_{ij}(\theta)=\begin{bmatrix}
H_{ij1}(\theta) & H_{ij2}(\theta)
\end{bmatrix}^{T}.
\end{align*}
Now, on the center manifold $C_{0}$, near the origin
 \begin{equation}
  \dot{\mathbf{w}} = \mathbf{w}_{z}\dot{z} + \mathbf{w}_{\bar{z}}\dot{\bar{z}}.\label{eq:derivative1_w}
 \end{equation} 
 We now use equations \eqref{eq:vector_w} and \eqref{eq:derivative_z} to replace $\mathbf{w}_{z}$ and $\dot{z}$ (and their conjugates) and equate this with \eqref{eq:derivative1_w} to get 
 \begin{align}
  (2i\omega_{0}- \mathcal{A})\mathbf{w}_{20}(\theta) &= \mathbf{H}_{20}(\theta),\notag\\
 -\mathcal{A}\mathbf{w}_{11}(\theta) &= \mathbf{H}_{11}(\theta),\notag\\
 -(2i\omega_{0}+ \mathcal{A})\mathbf{w}_{02}(\theta) &= \mathbf{H}_{02}(\theta)\label{eq:H_operator},
 \end{align}
as in \cite{Hassard}. Now, we observe that
 \begin{align}
\mathbf{u}_{t}(\theta) =&\,\, \mathbf{w}(z,\bar{z},\theta)+ z\mathbf{q}(\theta) + \bar{z}\bar{\mathbf{q}}(\theta)\notag\\
 =&\,\, \mathbf{w}_{20}(\theta)\frac{z^{2}}{2} + \mathbf{w}_{11}(\theta)z\bar{z} + \mathbf{w}_{02}(\theta)\frac{\bar{z}^{2}}{2}+ ze^{i\omega_{0}\theta} + \bar{z}e^{-i\omega_{0}\theta}+ \cdots, \label{eq:vector_u}
 \end{align}
 from which we obtain $\mathbf{u}_{t}(0)$, $\mathbf{u}_{t}(-\tau_1)$, and $\mathbf{u}_{t}(-\tau_2)$. We now proceed to expand the non-linear terms present in equation \eqref{eq:TaylorSeries} using equation \eqref{eq:vector_u} and retain only the coefficients of $z^2,z\bar{z},\bar{z}^2,z^2\bar{z}$. They are summarised as below:
 \begin{align*}
& u_{1,t}^2(0) = z^2 +\bar{z}^2 + 2z\bar{z} +z^2\bar{z}\Big( w_{201}(0) + 2w_{111}(0)\Big) + \cdots,\notag\\
& u_{2,t}^2(0) = \phi_1^2 z^2 +\bar{\phi_1}^{2}\bar{z}^2 + 2\phi_1 \bar{\phi_1} z\bar{z} + z^2\bar{z}\Big( w_{202}(0)\bar{\phi_1} + 2w_{112}(0)\phi_1\Big) + \cdots,\notag\\
& u_{1,t}^2(-\tau_1 ) =  z^2e^{-2i\omega_0 \tau_1} +\bar{z}^2 e^{2i\omega_0 \tau_1} + 2z\bar{z} 
 +z^2\bar{z}\Big( w_{201}(-\tau_1)e^{i\omega_0\tau_1} + 2w_{111}(-\tau_1)e^{-i\omega\tau_1} \Big) + \cdots,\notag\\
& u_{2,t}^2(-\tau_2 ) = \phi_1^2 z^2e^{-2i\omega_0 \tau_2} +\bar{\phi_1}^{2}\bar{z}^2 e^{2i\omega_0 \tau_2}+ 2\phi_1 \bar{\phi_1} z\bar{z}
 +z^2\bar{z}\Big( w_{202}(-\tau_2)\bar{\phi_1}e^{i\omega_0\tau_2}
 + 2w_{112}(-\tau_2)\phi_1e^{-i\omega\tau_2} \Big) + \cdots,\notag\\
& u_{1,t}(0)u_{1,t}(-\tau_1 ) =  z^2e^{-i\omega_0 \tau_1} +\bar{z}^2 e^{i\omega_0 \tau_1}
 + \Big(e^{i\omega_0 \tau_1}+e^{-i\omega_0 \tau_1}\Big)z\bar{z} 
 + z^2\bar{z}\Big(\frac{w_{201}(0)}{2}e^{i\omega_0\tau_1}
 +w_{111}(0)e^{-i\omega_0\tau_1}\notag\\
& \hspace*{25mm} +w_{111}(-\tau_1)+\frac{w_{201}(-\tau_1)}{2}\Big)\cdots,\notag\\
& u_{1,t}(0)u_{2,t}(-\tau_2 ) =  \phi_1 z^2e^{-i\omega_0 \tau_2} +\bar{\phi_1}\bar{z}^2 e^{i\omega_0 \tau_2}
 + \Big(\bar{\phi_1}e^{i\omega_0 \tau_2}+\phi_1e^{-i\omega_0 \tau_2}\Big)z\bar{z}
 + z^2\bar{z}\Big(\bar{\phi_1}\frac{w_{201}(0)}{2}e^{i\omega_0\tau_2} \notag\\
& \hspace*{25mm} +\phi_1 w_{111}(0)e^{-i\omega_0\tau_2}
 +w_{112}(-\tau_2)+\frac{w_{202}(-\tau_2)}{2}\Big)\cdots,\notag\\
& u_{1,t}(-\tau_1)u_{2,t}(-\tau_2 ) =  \phi_1 z^2e^{-i\omega_0(\tau_1+ \tau_2)} +\bar{\phi_1}\bar{z}^2 e^{i\omega_0 (\tau_1+\tau_2)}
 + \Big(\bar{\phi_1}e^{i\omega_0 (\tau_2-\tau_1)}+\phi_1e^{-i\omega_0 (\tau_2-\tau_1)}\Big)z\bar{z} \notag\\
& \hspace*{29mm} + z^2\bar{z}\Big(\bar{\phi_1}\frac{w_{201}(-\tau_1)}{2}e^{i\omega_0\tau_2}
 +\phi_1 w_{111}(-\tau_1)e^{-i\omega_0\tau_2}
+w_{112}(-\tau_2)e^{-i\omega_0\tau_1}
 +\frac{w_{202}(-\tau_2)}{2}e^{i\omega_0\tau_1}\Big)\cdots,\notag\\
  \end{align*}
 \begin{align*}
& u_{2,t}(0)u_{2,t}(-\tau_2 ) =  \phi_1^2 z^2e^{-i\omega_0 \tau_2} +\bar{\phi_1}^2\bar{z}^2 e^{i\omega_0 \tau_2}
 + \phi_1\bar{\phi_1}\Big(e^{i\omega_0 \tau_2}+e^{-i\omega_0 \tau_2}\Big)z\bar{z}
 + z^2\bar{z}\Big(\bar{\phi_1}\frac{w_{202}(0)}{2}e^{i\omega_0\tau_2} \notag \\
& \hspace*{24mm} +\phi_1 w_{112}(0)e^{-i\omega_0\tau_2} 
  +\phi_1 w_{112}(-\tau_2)  + \bar{\phi_1}\frac{w_{202}(-\tau_2)}{2}\Big)\cdots,\notag\\
& u_{2,t}(0)u_{1,t}(-\tau_1 ) =  \phi_1 z^2e^{-i\omega_0 \tau_1} +\bar{\phi_1}\bar{z}^2 e^{i\omega_0 \tau_1}
 + \Big(\phi_1 e^{i\omega_0 \tau_1}+\bar{\phi_1}e^{-i\omega_0 \tau_1}\Big)z\bar{z} + z^2\bar{z}\Big(\frac{w_{202}(0)}{2}e^{i\omega_0\tau_1} \notag\\
& \hspace*{24mm} + w_{112}(0)e^{-i\omega_0\tau_1}
 +\phi_1 w_{111}(-\tau_1)+\bar{\phi_1}\frac{w_{201}(-\tau_1)}{2}\Big)\cdots,\notag\\
&u_{1,t}^3(0) = 3z^{2}\bar{z}+ \cdots,\notag\\
& u_{2,t}^3(0) = 3\phi_1^2\bar{\phi_1} z^{2}\bar{z}+ \cdots,\notag\\
& u_{1,t}^3(-\tau_1) = 3z^{2}\bar{z}e^{-i\omega \tau_1}+ \cdots,\notag\\
& u_{2,t}^3(-\tau_2) = 3\phi_1^2\bar{\phi_1}z^{2}\bar{z}e^{-i\omega \tau_2}+ \cdots,\notag\\
& u_{1,t}^2(0)u_{1,t}(-\tau_1) =z^2\bar{z} \Big(2e^{-i\omega_0\tau_1}+e^{i\omega_0 \tau_1}\Big)\cdots,\notag\\
& u_{1,t}^2(0)u_{2,t}(-\tau_2) =z^2\bar{z} \Big(2\phi_1 e^{-i\omega_0\tau_2}+\bar{\phi_1}e^{i\omega_0 \tau_2}\Big)\cdots,\notag\\
& u_{1,t}^2(-\tau_1)u_{1,t}(0) =z^2\bar{z} \Big(e^{-2i\omega_0\tau_1}+2\Big)\cdots,\notag\\
& u_{1,t}^2(-\tau_1)u_{2,t}(-\tau_2) =z^2\bar{z} \Big(\bar{\phi_1}e^{-i\omega_0(2\tau_1-\tau_2)} +2\phi_{1}e^{-i\omega_0\tau_2}\Big)\cdots,\notag\\
 &  u_{2,t}^2(-\tau_2)u_{1,t}(0) =z^2\bar{z} \phi_1 \Big(\phi_1 e^{-2i\omega_0\tau_1}+2\bar{\phi_1}\Big)\cdots,\notag\\
&   u_{2,t}^2(-\tau_2)u_{1,t}(-\tau_1) =z^2\bar{z}\phi_1 \Big(\phi_1e^{-i\omega_0(2\tau_2-\tau_1)}+2\bar{\phi_1}e^{-i\omega_0\tau_1}\Big)\cdots,\notag\\
 &  u_{2,t}^2(0)u_{1,t}(-\tau_1) =z^2\bar{z} \phi_1\Big(2\bar{\phi_1}e^{-i\omega_0 \tau_1}+\phi_1e^{i\omega_0\tau_1}\Big)\cdots,\notag\\
& u_{2,t}^2(0)u_{2,t}(-\tau_2) =z^2\bar{z} \phi_1^2\bar{\phi_1}\Big(2e^{-i\omega_0 \tau_2}+e^{i\omega_0\tau_2}\Big)\cdots,\notag\\
& u_{2,t}^2(-\tau_2)u_{1,t}(0) =z^2\bar{z} \phi_1\Big(\phi_1e^{-2i\omega_0 \tau_2}+2\bar{\phi_1}\Big)\cdots,\notag\\
& u_{1,t}^2(-\tau_1)u_{2,t}(0) =z^2\bar{z} \Big(\bar{\phi_1}e^{-2i\omega_0 \tau_1}+2\phi_1\Big)\cdots,\notag\\
& u_{2,t}^2(-\tau_2)u_{2,t}(0) =z^2\bar{z}\phi_1^2\bar{\phi_1} \Big(2+e^{-2i\omega_0 \tau_2}\Big)\cdots,\notag\\
& u_{1,t}(0)u_{1,t}(-\tau_1)u_{2,t}(-\tau_2) =z^2\bar{z}\Big(\bar{\phi_1}e^{-i\omega_0 (\tau_1-\tau_2)} +\phi_1 e^{i\omega_0 (\tau_1-\tau_2)} +\phi_ 1e^{-i\omega_0 (\tau_1+\tau_2)}\Big)\cdots,\notag\\
&  u_{2,t}(0)u_{1,t}(-\tau_1)u_{2,t}(-\tau_2) =z^2\bar{z}\phi_1^2\Big(e^{-i\omega_0 (\tau_1-\tau_2)} + e^{i\omega_0 (\tau_1-\tau_2)} +e^{-i\omega_0 (\tau_1+\tau_2)}\Big)\cdots.
\end{align*}
Using the definition $g(z,\bar{z})=\bar{\mathbf{p}}(0)\cdot \mathcal{F}_{0}(z,\bar{z})$ we then determine the coefficients of $z^2$, $z\bar{z}$, $\bar{z}^2$ and $z^2\bar{z}$, which are outlined below.
\begin{align*}
 g_{20} &=2\kappa \bar{D} \Bigg( \bar{\phi}_{2}\Big(\xi_{aa}+\xi_{bb}e^{-2i\omega_0 \tau_1}+\xi_{dd}\phi_1^2 e^{-2i\omega_0 \tau_2}+ \xi_{ab}e^{-i\omega_0 \tau_1}+\xi_{ad}\phi_1 e^{-i\omega_0 \tau_2}+ \xi_{bd}\phi_1 e^{-i\omega_0(\tau_1+ \tau_2)}\Big)\notag\\
 &+ \chi_{cc}\phi_1^2 +\chi_{dd}\phi_1^2 e^{-2i\omega_0 \tau_2}+\chi_{bb}e^{-2i\omega_0 \tau_1}+ \chi_{cd}\phi_1^2 e^{-i\omega_0 \tau_2} +\chi_{bc}\phi_1 e^{-i\omega_0 \tau_1}+ \chi_{bd}\phi_1 e^{-i\omega_0(\tau_1+ \tau_2)}\Bigg),\notag\\
 \end{align*}
 \begin{align}
 g_{02} &=2\kappa \bar{D} \Bigg( \bar{\phi}_{2}\Big(\xi_{aa}+\xi_{bb}e^{2i\omega_0 \tau_1}+\xi_{dd}\phi_1^2 e^{2i\omega_0 \tau_2}+ \xi_{ab}e^{i\omega_0 \tau_1} +\xi_{ad}\bar{\phi_1} e^{i\omega_0 \tau_2}+ \xi_{bd}\bar{\phi_1} e^{i\omega_0(\tau_1+ \tau_2)}\Big)\notag\\
 &+ \chi_{cc}\bar{\phi_1}^2 +\chi_{dd}\bar{\phi_1}^2 e^{2i\omega_0 \tau_2}+\chi_{bb}e^{2i\omega_0 \tau_1}+ \chi_{cd}\bar{\phi_1}^2 e^{i\omega_0 \tau_2} +\chi_{bc}\bar{\phi_1} e^{i\omega_0 \tau_1}+ \chi_{bd}\bar{\phi_1} e^{i\omega_0(\tau_1+ \tau_2)}\Bigg),\notag\\
  g_{11} &=\kappa \bar{D} \Bigg( \bar{\phi}_{2}\Big(2\xi_{aa}+2\xi_{bb}+\xi_{dd}\phi_1\bar{\phi_1}+ \xi_{ab}\left(e^{i\omega_0 \tau_1}+e^{-i\omega_0 \tau_1}\right)+\xi_{ad}\left(\bar{\phi_1}e^{i\omega_0 \tau_2}+\phi_1e^{-i\omega_0 \tau_2}\right)\notag\\
 &+ \xi_{bd}\left(\bar{\phi_1}e^{i\omega_0 (\tau_2-\tau_1)}+\phi_1e^{-i\omega_0 (\tau_2-\tau_1)}\right)\Big)+2\phi_1\bar{\phi_1}(\chi_{cc} +\chi_{dd}) +2\chi_{bb}+\chi_{cd}\phi_1\bar{\phi_1}\left(e^{i\omega_0 \tau_2}+e^{-i\omega_0 \tau_2}\right)\notag\\ &+\chi_{bc}\left(\phi_1 e^{i\omega_0 \tau_1}+\bar{\phi_1}e^{-i\omega_0 \tau_1}\right)+ \chi_{bd}\left(\bar{\phi_1}e^{i\omega_0 (\tau_2-\tau_1)}+\phi_1e^{-i\omega_0 (\tau_2-\tau_1)}\right)\Bigg),\notag\\
g_{21}&=2\kappa \bar{D}\Bigg(\bar{\phi_2} \Big(\xi_{aa}\left(w_{201}(0) + 2w_{111}(0)\right)+ \xi_{bb}\left(w_{201}(-\tau_1)e^{i\omega_0\tau_1}+2w_{111}(-\tau_1)e^{-i\omega\tau_1}\right)\notag\\
&+ \xi_{dd}\left(w_{202}(-\tau_2)\bar{\phi_1}e^{i\omega_0\tau_2} +2w_{112}\phi_1(-\tau_2)e^{-i\omega\tau_2}\right)+ \xi_{ab}\big(\frac{w_{201}(0)}{2}e^{i\omega_0\tau_1}+w_{111}(0)e^{-i\omega_0\tau_1}+ w_{111}(-\tau_1)\notag\\
&+\frac{w_{201}(-\tau_1)}{2}\big)+ \xi_{ab}\big(\frac{w_{201}(0)}{2}\bar{\phi_1}e^{i\omega_0\tau_2}+w_{111}(0)\phi_1e^{-i\omega_0\tau_2}+ w_{112}(-\tau_2)+\frac{w_{202}(-\tau_2)}{2}\big)\notag\\
& + \xi_{bd}\big(\frac{w_{201}(-\tau_1)}{2}\bar{\phi_1}e^{i\omega_0\tau_2}+ w_{111}(-\tau_1)\phi_1e^{-i\omega_0\tau_2}+ w_{112}(-\tau_2)e^{-i\omega_0\tau_1}+\frac{w_{202}(-\tau_2)}{2}e^{i\omega_0\tau_1}\big)\notag\\
&+3\xi_{aaa}+3\xi_{bbb}e^{i\omega_0\tau_1}
+3\xi_{ddd}\phi_1^2\bar{\phi_1}e^{-\i\omega_0\tau_2}+\xi_{aab} \big(2e^{-i\omega_0\tau_1}+e^{i\omega_0 \tau_1}\big)+\xi_{aad}\big(2\phi_1 e^{-i\omega_0\tau_2}+\bar{\phi_1}e^{i\omega_0 \tau_2}\big)\notag\\
&+\xi_{abb} \big(e^{-2i\omega_0\tau_1}+2\big)+\xi_{bbd}\big(\bar{\phi_1}e^{-i\omega_0(2\tau_1-\tau_2)}+2\phi_{1}e^{-i\omega_0\tau_2}\big)+ \xi_{add}\phi_1 \big(\phi_1e^{-2i\omega_0\tau_2}+2\bar{\phi_1}\big)\notag\\
 &+\xi_{bdd}\phi_1 \big(\phi_1e^{-i\omega_0(2\tau_2-\tau_1)}+\bar{\phi_1}e^{-i\omega_0\tau_1}\big)+\xi_{abd}\big(\bar{\phi_1}e^{-i\omega_0 (\tau_1-\tau_2)}+\phi_1 e^{i\omega_0 (\tau_1-\tau_2)}+\phi_ 1e^{-i\omega_0 (\tau_1+\tau_2)}\big)\Big)\notag\\
 &+ \chi_{cc}\left(w_{202}(0)\bar{\phi_1} + 2w_{112}(0)\phi_1\right)+\chi_{dd}\left(w_{202}(-\tau_2)\bar{\phi_1}e^{i\omega_0\tau_2}+2w_{112}(-\tau_2)\phi_1e^{-i\omega\tau_2}\right)\notag\\
&+ \chi_{bb}\left(w_{201}(-\tau_1)e^{i\omega_0\tau_1}+2w_{111}(-\tau_1)e^{-i\omega\tau_1}\right)+ \chi_{cd}\big(\frac{w_{202}(0)}{2}\bar{\phi_1}e^{i\omega_0\tau_2}+\phi_1 w_{112}(0)e^{-i\omega_0\tau_2}\notag\\
&+\phi_1 w_{112}(-\tau_2)+\bar{\phi_1}\frac{w_{202}(-\tau_2)}{2}\big)+\chi_{bc}\big(\frac{w_{202}(0)}{2}e^{i\omega_0\tau_1}+ w_{112}(0)e^{-i\omega_0\tau_1}+\phi_1 w_{111}(-\tau_1)+\bar{\phi_1}\frac{w_{201}(-\tau_1)}{2}\big)\notag\\
 &+\chi_{bd}\big(\bar{\phi_1}\frac{w_{201}(-\tau_1)}{2}e^{i\omega_0\tau_2}+\phi_1 w_{111}(-\tau_1)e^{-i\omega_0\tau_2}+w_{112}(-\tau_2)e^{-i\omega_0\tau_1}+\frac{w_{202}(-\tau_2)}{2}e^{i\omega_0\tau_1}\big)+3\chi_{ccc}\phi_1^2\bar{\phi_1}\notag\\ &+3\chi_{ddd}\phi_1^2\bar{\phi_1}e^{-i\omega_0\tau_2}+3\chi_{bbb}e^{-i\omega_0\tau_1}+ \chi_{ccd}\phi_1^2\bar{\phi_1}\big(2e^{-i\omega_0 \tau_2}+e^{i\omega_0\tau_2}\big)+\chi_{bcc}\phi_1 \big(2\bar{\phi_1}e^{-i\omega_0-\tau_1}+\phi_1e^{i\omega_0\tau_1}\big)\notag\\
 &+\chi_{cdd} \phi_1^2\bar{\phi_1}\big(e^{-2i\omega_0 \tau_2}+2\big)+\chi_{bdd}\phi_1 \big(\phi_1e^{-i\omega_0(2\tau_2-\tau_1)}+2\bar{\phi_1}e^{-i\omega_0\tau_1}\big)+ \chi_{bbc}\big(\bar{\phi_1}e^{-2i\omega_0 \tau_2}+2\phi_1\big)\notag\\
 &+\chi_{bbd}\big(\bar{\phi_1}e^{-i\omega_0(2\tau_1-\tau_2)}+2\phi_{1}e^{-i\omega_0\tau_2}\big)+\chi_{bcd}\phi_1\big(\bar{\phi_1}e^{-i\omega_0 (\tau_1-\tau_2)}+\phi_1 e^{i\omega_0 (\tau_1-\tau_2)}+\bar{\phi_1}e^{-i\omega_0 (\tau_1+\tau_2)}\big)\Bigg). 
\label{eq:g21}
\end{align}
Note that, the expression for $g_{21}$ has $\mathbf{w}_{20}(\theta)$ and $\mathbf{w}_{11}(\theta)$ which we need to evaluate. Now, for $\theta \in [-\tau,0)$ from \eqref{eq:vector_H}, we have
\begin{align*}
   \mathbf{H}(z,\bar{z},\theta) =& -2\text{Re}\big(\mathbf{\bar{q}}^{\ast}(0)\cdot \mathcal{F}_{0} \mathbf{q}(\theta)\big)\\
 =& -g(z,\bar{z})\mathbf{q}(\theta) - \bar{g}(z,\bar{z})\mathbf{\bar{q}}(\theta)\\
   =& -\left(g_{20}\frac{z^{2}}{2} + g_{11}z\bar{z} + g_{02}\frac{\bar{z}^{2}}{2}+\cdots \right)\mathbf{q}(\theta) - \left(\bar{g}_{20}\frac{\bar{z}^{2}}{2} + \bar{g}_{11}z\bar{z} + \bar{g}_{02}\frac{z^{2}}{2}+\cdots \right)\mathbf{\bar{q}}(\theta),
  \end{align*}
which when compared with \eqref{eq:vector_H} gives
 \begin{align}
 \label{eq:vector1_H}
   \mathbf{H}_{20}(\theta) &= -g_{20}\mathbf{q}(\theta)-\bar{g}_{02}\mathbf{\bar{q}}(\theta),\notag\\
   \mathbf{H}_{11}(\theta) &= -g_{11}\mathbf{q}(\theta)-\bar{g}_{11}\mathbf{\bar{q}}(\theta).
   \end{align}
Using equations \eqref{eq:DLEoperatorR} and \eqref{eq:H_operator}, we have
 \begin{align}\label{eq:w_derivatives}
  \dot{\mathbf{w}}_{20}(\theta) &= 2i\omega_{0}\mathbf{w}_{20}(\theta) + g_{20}\mathbf{q}(\theta)+ \bar{g}_{02}\mathbf{\bar{q}}(\theta),\notag\\
   \dot{\mathbf{w}}_{11}(\theta) &= g_{11}\mathbf{q}(\theta) + \bar{g}_{11}\mathbf{\bar{q}}(\theta).
  \end{align}
  Solving the differential equations in \eqref{eq:w_derivatives}, we get
 \begin{align}\label{eq:w_vectorssol}
  \mathbf{w}_{20}(\theta) &= -\frac{g_{20}}{i\omega_{0}}\boldsymbol{q}(0)e^{i\omega_{0}\theta}-\frac{\bar{g}_{02}}{3i\omega_{0}}\mathbf{\bar{q}}(0)e^{-i\omega_{0}\theta} + \mathbf{e}e^{2i\omega_{0}\theta},\notag\\
  \mathbf{w}_{11}(\theta) &= \frac{g_{11}}{i\omega_{0}}\mathbf{q}(0)e^{i\omega_{0}\theta} - \frac{\bar{g}_{11}}{i\omega_{0}}\mathbf{\bar{q}}(0)e^{-i\omega_{0}\theta} + \boldsymbol{f}.
  \end{align}
The objective now is to determine $ \boldsymbol{e}$ and $ \boldsymbol{f}$. We define,
\begin{align}
  \mathbf{H}(z,\bar{z},0) &= -2\text{Re}\big( \mathbf{\bar{q}}^{\ast}(0)\cdot\mathcal{F}_{0} \mathbf{q}(0)\big)+\mathcal{F}_0,\label{eq:H_0}
 \end{align}
where $\boldsymbol{\mathcal{{F}}}_0$ represents the non-linear terms that can be expanded in powers of $z$ as
  \begin{equation}
   \boldsymbol{\mathcal{F}}=\boldsymbol{\mathcal{F}}_{20}\frac{z^2}{2}+ \boldsymbol{\mathcal{F}}_{11}z\bar{z}+\boldsymbol{\mathcal{F}}_{02}\frac{\bar{z}^2}{2} + \boldsymbol{\mathcal{F}}_{21}\frac{z^2 \bar{z}}{2} + \cdots.
  \end{equation}
  Substituting the coefficients from the expansion of $\boldsymbol{\mathcal{F}}_0$ gives
 \begin{align}
 \label{eq:H_coeff}
 \mathbf{H}_{20}(0) &= -g_{20}\mathbf{q}(0)-\bar{g}_{02}\mathbf{\bar{q}}(0)+\begin{bmatrix}\mathcal{F}_{201}&\mathcal{F}_{202}\end{bmatrix}^{T},\notag\\
\mathbf{H}_{11}(0) &= -g_{11}\mathbf{q}(0)-\bar{g}_{11}\mathbf{\bar{q}}(0)+\begin{bmatrix}\mathcal{F}_{111}&\mathcal{F}_{112}\end{bmatrix}^{T}.
\end{align}
From \eqref{eq:H_coeff} and \eqref{eq:DLEoperatorR}, we obtain
\begin{align}
 &g_{20}\mathbf{q}(0)+\bar{g}_{02}\mathbf{\bar{q}}(0) = \begin{bmatrix}\mathcal{F}_{201}&\mathcal{F}_{202}\end{bmatrix}^{T}+ \begin{bmatrix} (\kappa a_{11}-2i\omega_0)\mathsf{w}_{201}(0) + \kappa a_{12}\mathsf{w}_{201}(-\tau_1)+\kappa a_{13}\mathsf{w}_{202(-\tau_1)}\notag\\
 \kappa a_{23}\mathsf{w}_{201}(-\tau)+\kappa(a_{21}- 2i\omega_0)\mathsf{w}_{202}(0)+\kappa a_{22}\mathsf{w}_{202}(-\tau_2)\end{bmatrix},\notag\\
&g_{11}\mathbf{q}(0)+\bar{g}_{11}\mathbf{\bar{q}}(0) = \begin{bmatrix}\mathcal{F}_{111}&\mathcal{F}_{112}\end{bmatrix}^{T}+ \begin{bmatrix}\kappa a_{11}\mathsf{w}_{111}(0)-\kappa a_{12}\mathsf{w}_{111}(-\tau_1)+\kappa a_{13}\mathsf{w}_{112}(-\tau_2)\\ \kappa a_{23}\mathsf{w}_{111}(-\tau_1)+\kappa a_{21} \mathsf{w}_{112}(0)+\kappa a_{22}\mathsf{w}_{112}(-\tau_2)\end{bmatrix}.\label{eq:getE&F}
\end{align}
We substitute $\mathbf{w}_{20}(0),\mathbf{w}_{20}(-\tau),\mathbf{w}_{11}(0)$ and $\mathbf{w}_{11}(-\tau)$ from \eqref{eq:w_vectorssol} in \eqref{eq:getE&F} we get $\mathbf{e}$ and $\mathbf{f}$ of the form
\begin{align}
 \boldsymbol{e} = \begin{bmatrix}e_{1}& e_{2}\end{bmatrix}^{T} \quad\text{and}\quad \boldsymbol{f} = \begin{bmatrix}f_{1}& f_{2}\end{bmatrix}^{T}.
\end{align}
Note that, $e_1,e_2,f_1$ and $f_2$ can be derived explicitly in terms of system parameters, which are outlined below:
\begin{align}
 e_1 = \frac{Y_2Z_1 - Y_1Z_2}{X_1Y_2-X_2Y_1}, \hspace{3mm} e_2 =\frac{X_1Z_2 - X_2Z_1}{X_1Y_2-X_2Y_1},\hspace{3mm} f_1 = \frac{Q_2R_1 - Q_1R_2}{P_1Q_2-P_2Q_1}, \hspace{3mm} f_2 = \frac{P_1R_2-P_2R_1}{P_1Q_2-P_2Q_1},\notag\\
\end{align}
where,
\begin{align}
X_1 &= \,\,  \kappa a_{11}+\kappa a_{12}e^{-2i\omega_0\tau_1}-2i\omega_0,\hspace{5mm} X_2 = \kappa a_{23}e^{-2i\omega _{0}\tau_1},\notag\\
Y_1 &= \,\,\kappa a_{13}e^{-2i\omega_0\tau_2}, \hspace{4mm}  Y_2 =  \kappa a_{21}+\kappa a_{22}e^{-2i\omega_0\tau_2}-2i\omega_0,\notag\\
P_1 &=  \,\, \kappa a_{11}+\kappa a_{12}, \hspace{5mm} P_2 = \kappa a_{23},\hspace{5mm}Q_1 =  \,\, \kappa a_{13}, \hspace{5mm} Q_2 = \kappa a_{21}+\kappa a_{22},\notag\\
Z_1 &= \,\, \frac{g_{20}}{i\omega_0}\left(-i\omega_0 + \kappa a_{11}+ \kappa a_{12}e^{-i\omega_0\tau_1}+\kappa a_{13}\phi_1 e^{-i\omega_0\tau_1} \right)+\frac{\bar{g}_{02}}{3i\omega_0}\left(i\omega_0 + \kappa a_{11}+ \kappa a_{12}e^{i\omega_0\tau_1}+\kappa a_{13}\bar{\phi}_1 e^{i\omega_0\tau_1}\right)\notag\\
&- \mathcal{F}_{201},\notag\\
Z_2 &= \,\, \frac{g_{20}}{i\omega_0}\big(-i\omega_0 \phi_1 + \kappa a_{23}e^{i \omega_{0}\tau_1}+\kappa a_{21}\phi_1+ \kappa a_{22}\phi_1 e^{i\omega_0\tau_2}\big)+\frac{\bar{g}_{02}}{3i\omega_0}\big(i\omega_0\bar{\phi}_1 + \kappa a_{23}e^{i \omega_{0}\tau_1}+\kappa a_{21}\bar{\phi}_1\notag\\
 &+ \kappa a_{22}\bar{\phi}_1 e^{i\omega_0\tau_2}\big)-\mathcal{F}_{202},\notag
\end{align}
\begin{align}
 R_1 &= \,\, \frac{g_{11}}{i\omega_0}\left(i\omega_0 -\kappa a_{11}- \kappa a_{12}e^{-i\omega_0\tau_1}-\kappa a_{13}\phi_1 e^{-i\omega_0\tau_1}\right)+ \frac{\bar{g}_{11}}{i\omega_0}\left(i\omega_0 + \kappa a_{11}+ \kappa a_{12}e^{i\omega_0\tau_1}+\kappa a_{13}\bar{\phi}_1 e^{i\omega_0\tau_1}\right)\notag\\
&-\mathcal{F}_{111},\notag\\
 R_2 &= \,\, \frac{g_{11}}{i\omega_0}\big(i\omega_0 \phi_1  \kappa a_{23}e^{i \omega_{0}\tau_1}-\kappa a_{21}\phi_1- \kappa a_{22}\phi_1 e^{i\omega_0\tau_2}\big)+\frac{\bar{g}_{02}}{3i\omega_0}\big(i\omega_0\bar{\phi}_1 + \kappa a_{23}e^{i \omega_{0}\tau_1}+\kappa a_{21}\bar{\phi}_1+ \kappa a_{22}\bar{\phi}_1 e^{i\omega_0\tau_2}\big)\notag\\
 &-\mathcal{F}_{112}.
 \end{align}
Using $\boldsymbol{e}$ and $\boldsymbol{f}$ we evaluate $\mathbf{w}_{20}$ and $\mathbf{w}_{11}$, using which we compute $g_{21}$. We now have all the terms required for the analysis of Hopf bifurcation as follows, see \cite{Hassard} 
\begin{align}
 \hspace{-3mm}c_1(0) &= \frac{i}{2\omega_0}\left(g_{20}g_{11}-2|g_{11}|^2-\frac{1}{3}|g_{02}|^2\right)+\frac{g_{21}}{2},\label{eq:cterm}\\
 \mu_2 &= -\frac{\text{Re}\big(c_1(0)\big)}{\alpha'(0)},\quad\quad\beta_2 = 2\text{Re}\big(c_1(0)\big),\label{eq:muterm&betaterm}
\end{align}
where $c_1(0)$ is the lyapunov coefficient and $\alpha'(0) =~\text{Re}\left(\mathrm{d}\lambda/\mathrm{d}\kappa\right)|_{\kappa=\kappa_c}$.
The following conditions enable us to verify the type of the Hopf bifurcation, and the asymptotic orbital stability of the limit cycles~\cite{Hassard}.
\begin{itemize}
 \item The Hopf bifurcation is \emph{supercritical} if $\mu_2 > 0$ and \emph{sub-critical} if $\mu_2 <0$.
 \item The limit cycles are \emph{asymptotically orbitally stable} if $\beta_2< 0$ and \emph{unstable} if $\beta_2>0$.
\end{itemize}
Substituting the expression for $g_{21}$ in \eqref{eq:cterm} yields the expression for $c_1(0)$, which is the lyapunov coefficient. We can then compute $\mu_2$ and $\beta_2$ using \eqref{eq:muterm&betaterm}. We now present a numerical example, and compute the values of $\mu_2$ and $\beta_2$ for Compound TCP in the small buffer regime.

\begin{figure*}
 \begin{center}
\begin{minipage}[b]{0.45\textwidth}
  \label{fig:phase_0.95}
  \psfrag{b}{$w_{2}(t)$}
  \psfrag{T}{$w_{2}(t-\tau_2)$}
    \psfrag{124}{\begin{scriptsize}$124$\end{scriptsize}}
  \psfrag{130}{\begin{scriptsize}$130$\end{scriptsize}}  
  \psfrag{136}{\begin{scriptsize}$136$\end{scriptsize}}
     \includegraphics[width=2in,height=3in,angle=270]{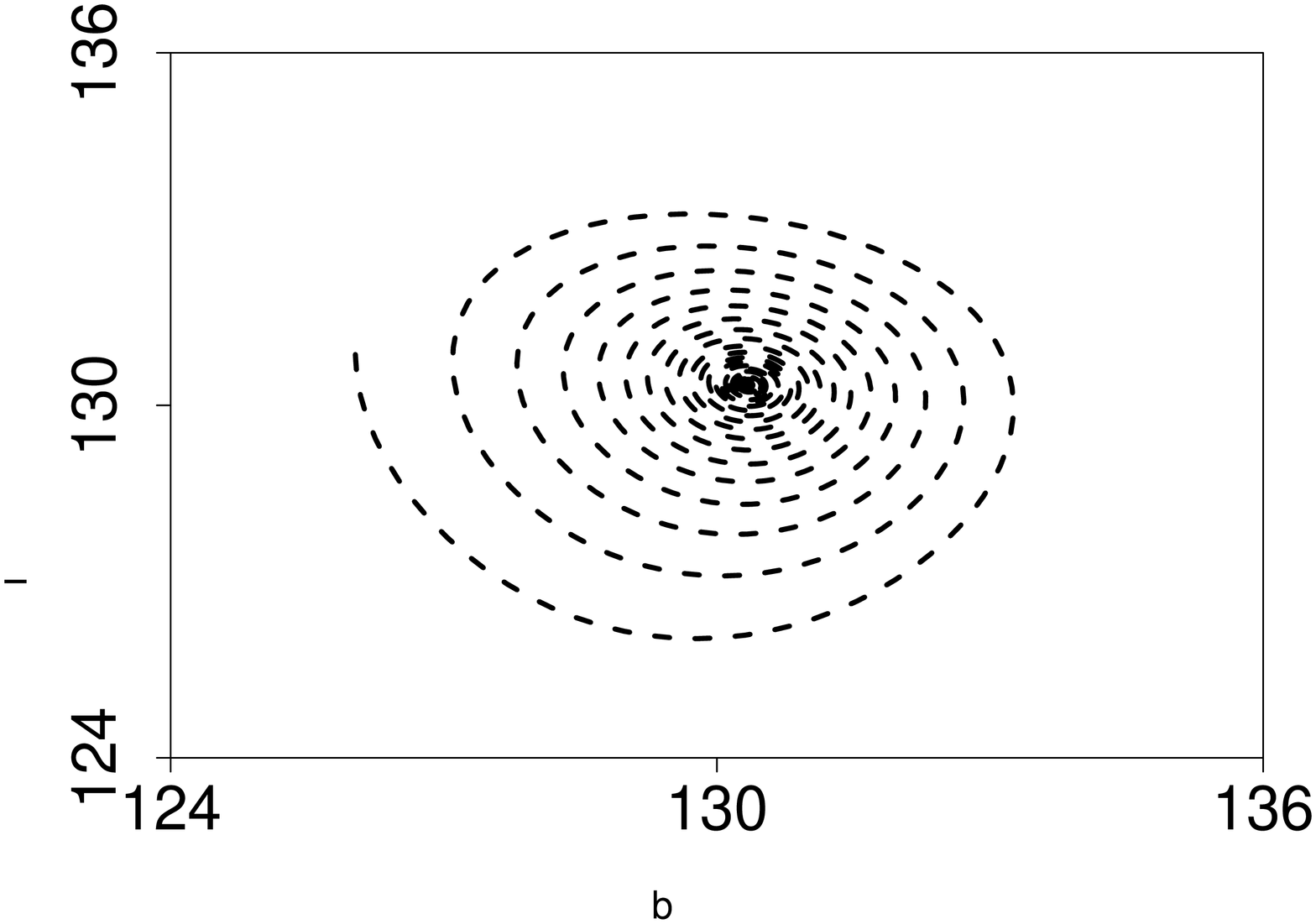}
  
\end{minipage}
\hspace*{5mm}
 \begin{minipage}[b]{0.45\textwidth}
   \label{fig:phase_1.05}
  \psfrag{b}{$w_2(t)$}
  \psfrag{T}{$w_2(t-\tau_2)$}
    \psfrag{122}{\begin{scriptsize}$122$\end{scriptsize}}
    \psfrag{129}{\begin{scriptsize}$129$\end{scriptsize}}
  \psfrag{136}{\begin{scriptsize}$136$\end{scriptsize}}

  \includegraphics[width=2in,height=3in,angle=270]{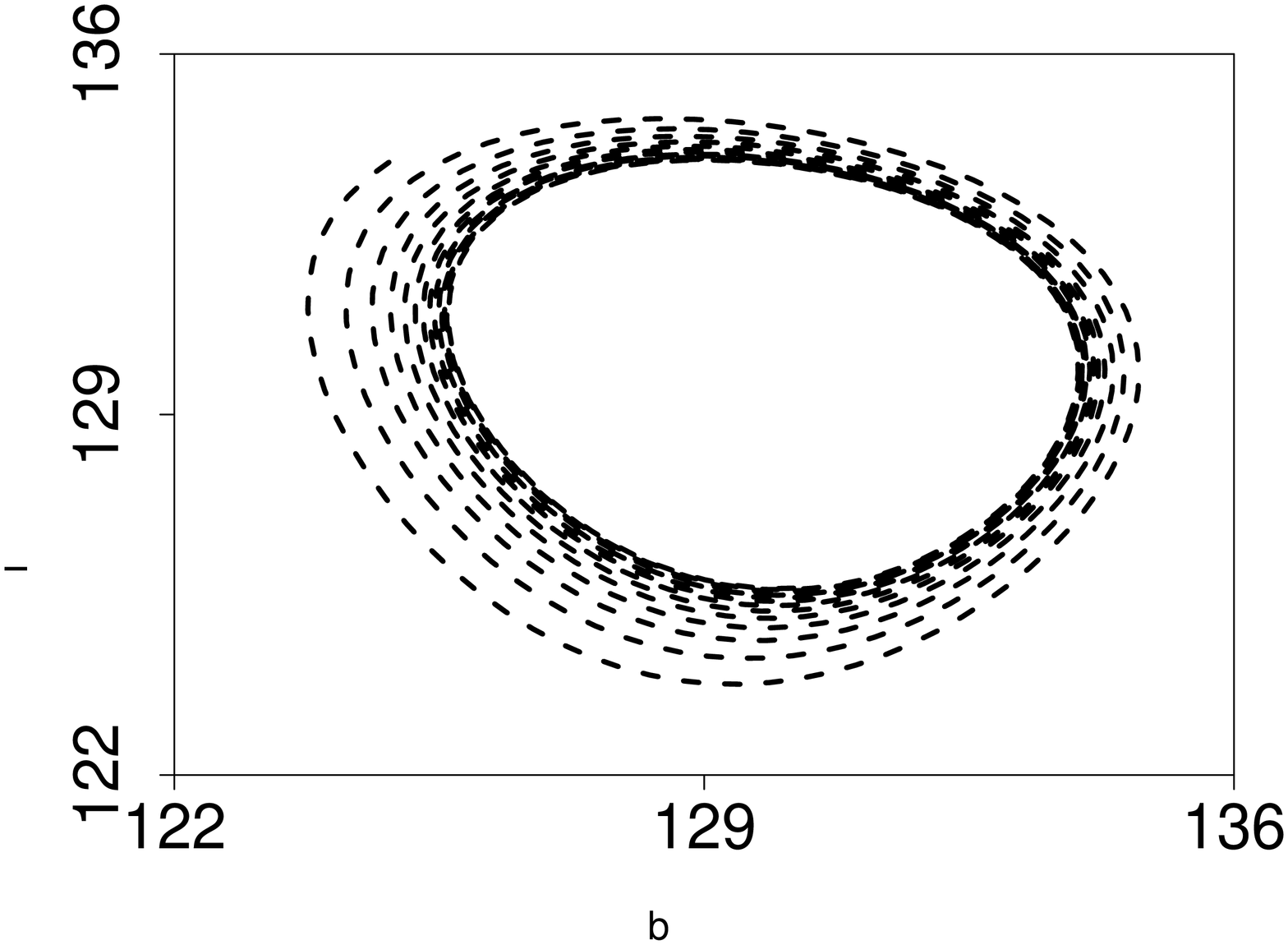}
  
\end{minipage}
  
     \caption {\emph{Phase portraits.} Emergence of limit cycle in the dynamics of $w_2(t)$ in \eqref{eq:modelb_3}, for Compound TCP in the small buffer regime, with the variation in the non-dimensional parameter $\kappa$. Observe that, (a) Trajectories converge to stable equilibrium for $\kappa=0.95$, (b) Trajectories converge to a stable limit cycle for $\kappa=1.05$.}
  \label{fig:phase}
 \end{center}
\end{figure*}

\subsection*{Numerical Example:}We first fix the system parameters as follows: $\alpha=0.3$, $k=0.75$, $\beta=0.5$, $B_1=10$, $B_2=15$, $B=25$, $C_1=C_2=100$, $C=180$, $\tau_1=1$, and $\tau_2$. With these parameter values, the system undergoes a Hopf bifurcation at $\kappa_c=1$. We now increase the value of the non-dimensional parameter to $\kappa=1.05$, and push the system beyond the edge of stability. Following the Hopf bifurcation analysis presented above, we compute the required expressions:
\begin{align*}
\mathrm{Re}\left(c_{1}(0)\right) &= -0.0738<0, \hspace{1ex} \alpha'(0)=0.3467>0\notag\\
\mu_2 &=0.2129>0, \hspace{5.5ex}\beta_{2}=-0.1477<0.
\end{align*}
Thus, the Hopf bifurcation is \emph{supercritical} and the emergent limit cycles are asymptotically \emph{orbitally stable}.
\subsection*{Phase portraits and bifurcation diagram:} 
\begin{figure}[t]
\begin{center}
 \psfrag{T}[0,0.5][0.25,0.75]{Amplitude}
  \psfrag{b}{\hspace{-21.5mm}Non-dimensional parameter, $\kappa$}
  \psfrag{0.9}{\begin{scriptsize}$0.9$\end{scriptsize}}
  \psfrag{1.1}{\begin{scriptsize}$1.1$\end{scriptsize}}
  \psfrag{1}{\begin{scriptsize}$1$\end{scriptsize}}
  \psfrag{120}{\begin{scriptsize}\hspace{1mm}$120$\end{scriptsize}}
  \psfrag{137}{\begin{scriptsize}$137$\end{scriptsize}}
      \psfrag{SSSSSSSSSSSSSS}{\hspace{2mm}\begin{scriptsize}Stable Region\end{scriptsize}}
  \psfrag{h}{\begin{scriptsize}Hopf Condition\end{scriptsize}}
  \includegraphics[width=1.8in,height=2.6in,angle=270]{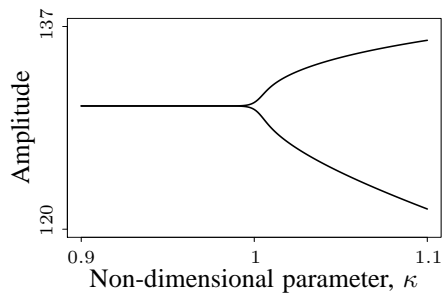}
  \caption{\emph{Bifurcation diagram.} Emergence of limit cycles in the dynamics of $w_2(t)$ at $\kappa=1$ for system \eqref{eq:modelb_3}, with Compound TCP flows in the small buffer regime. The amplitude of the emergent limit cycles increases for further increase in $\kappa$. }\label{fig:bifurcation}\vspace{-8mm}
  \end{center}
\end{figure}

We present the phase portrait for system \eqref{eq:modelb_3}, for Compound TCP in the small buffer regime, in Fig.~\ref{fig:phase}. First, we fix a point $\alpha=0.3,$ $\kappa=1$, on the stability boundary in the stability chart as shown in  Fig.~\ref{fig:charts}~(a). The remaining system parameter values are fixed as mentioned above in the numerical example. We now plot the phase portrait for the window size for the second set of TCP flows, for $\kappa=0.95$ and $\kappa=1.05$ respectively, as shown in Fig.~\ref{fig:phase}. Observe that, for $\kappa=0.95$, the average window size of the second set Compound TCP flows converges to its equilibrium value, as expected. For, $\kappa=1.05$, the average window size exhibits orbitally stable limit cycles, as the system undergoes a Hopf bifurcation at $\kappa=1$. Note that, the average window size of the first set of Compound flows can be shown to exhibit qualitatively similar dynamical behaviour. We now present the bifurcation diagram for system \eqref{eq:modelb_3}, in Fig.~\ref{fig:bifurcation}, obtained from DDE-BIFTOOL version $2.03$. Observe that, the amplitude of the limit cycles increases as $\kappa$ is increased beyond $1$.

\section{Packet-level simulations}
\label{simulations}
In order to corroborate the analytical insights obtained, we conduct some packet-level simulations, for the multiple bottleneck scenario, in NS2 \cite{ns2}.
\begin{figure}
\begin{center}
   \psfrag{0}{\begin{scriptsize}$0$\end{scriptsize}}
  \psfrag{15}{\begin{scriptsize}$15$\end{scriptsize}}
  \psfrag{100}{\begin{scriptsize}$100$\end{scriptsize}}
  \psfrag{125}{\begin{scriptsize}$125$\end{scriptsize}}
      \psfrag{aaaa}{\begin{scriptsize}Buffer size = 15 pkts\end{scriptsize}}
  \psfrag{bbbb}{\begin{scriptsize}Buffer size = 100 pkts\end{scriptsize}}
\psfrag{cccc}{Queue size (pkts)}
\psfrag{xyz}{\hspace{-6mm}Time (seconds)}
  \psfrag{mmmm}{\begin{scriptsize}$\tau_1=10$ ms, $\tau_2= 10$ ms\end{scriptsize}}
   \psfrag{nnnn}{\begin{scriptsize}\hspace{-1mm}$\tau_1=10$ ms, $\tau_2= 200$ ms\end{scriptsize}}

  \includegraphics[width=2.5in,height=3.5in,angle=270]{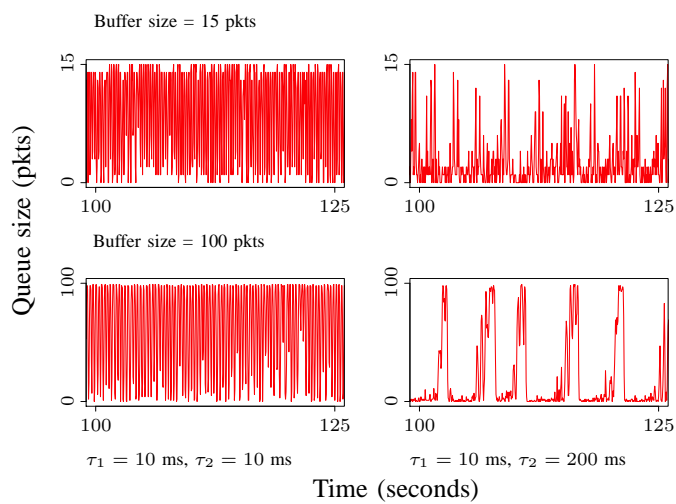}
  \caption{ \emph{Long-lived flows}. Two sets of 60 long-lived Compound flows over a 2 Mbps link, regulated by two edge routers, feeding into a core router with link capacity 180 Mbps. Observe the emergence of limit cycles in the queue at the core router, for larger buffer thresholds, and larger round trip times.}\vspace{-8mm}
  \label{ns2}
  \end{center}
\end{figure}
The system consists of two distinct sets of $60$ long-lived Compound TCP flows each with an access speed of $2$ Mbps, regulated by two edge routers and feeding into one core router. Each edge router has a link capacity of $100$ Mbps, and the core router has a link capacity of $180$ Mbps. Since our primary focus is on small buffers, we fix the buffer size for each edge router to be $15$ packets, and vary the buffer size of the core router from $15$ packets to $100$. Further, we fix the round trip time of one set of flows to be $10$ ms, and the round trip time of the other set is varied from $10$ ms to $200$ ms. The simulations are illustrated in Fig.~\ref{ns2}. Observe that, if the buffer sizes at all routers are fixed at $15$ packets, the queue at the core router is completely random, and hence stable, since the queue does not exhibit any deterministic oscillations. When the buffer size at the core router is increased to $100$ packets and the round trip time of the second set of flows is $200$ ms, the queue dynamics exhibits limit cycles. Hence, larger queue thresholds are prone to inducing limit cycles, for larger round trip times. These limit cycles in the queue size lead to synchronisation among TCP flows and make the downstream traffic bursty.  
\section{Concluding remarks}

We considered three different topologies, and conducted a detailed local stability analysis with two simplifying assumptions, to obtain necessary and sufficient conditions for stability. To aid our analysis, we motivated a suitable non-dimensional bifurcation parameter, and illustrated that, the underlying dynamical systems lose stability if the bifurcation parameter is varied. Further, in the multiple bottleneck scenario, even without any simplifying assumptions on the system parameters, we numerically identified that the system loses stability via a Hopf bifurcation. A key insight obtained was the trade-off between different system parameters to ensure stability, as illustrated through some stability charts. After knowing that a system exhibits a Hopf, it is natural to have a framework to determine the asymptotic orbital stability of the bifurcating limit cycles. To that end, using Poincar\'{e} normal forms and the center manifold theory, we conducted a detailed Hopf bifurcation analysis, in the neighbourhood of the Hopf condition. To corroborate our analytical insights, we conducted some packet-level simulations to highlight the existence and stability of limit cycles in the queue size dynamics as system parameters vary.          

The insights obtained in this paper could have important consequences for the modelling and the  performance evaluation of communication networks. From a theoretical perspective, this opens many challenging questions centred around the development of accurate fluid models for TCP and queue management policies. From a practical perspective, the emergence of stable limit cycles could have an impact on the end-to-end quality of service -- these issues merit further investigation.

\label{conclusions}

\begin{table*}\section*{Appendix}
\centering
\setlength{\tabcolsep}{36pt}
\caption{Coefficients in the Taylor series expansion of the non-linear fluid model~\eqref{eq:modelb_3} with Compound TCP and Drop-tail queue policy evaluated at the equilibrium $(w_1^\ast, w_2^\ast)$. Here, the term $p'$ represents the partial derivative of $p$ with respect to the variables as given by the subscripts.}
\begin{tabular}{l l }\hline\hline
\\ [-1mm]
 $\displaystyle\xi_a=\frac{\partial f_1}{\partial w_{1}(t)}$ & $\displaystyle\frac{w_1^\ast}{\tau_1}\left(\frac{ i_{1,a}^{'}d_{1}^{*}-i_{1}^{\ast}   d_{1,a}^{'}}{i_1^{\ast}+d_1^{\ast}}\right)$ \\[3ex]
 $\displaystyle\xi_b=\frac{\partial f_1}{\partial w_{1}(t-\tau_1)}$ & $\displaystyle -\frac{w_1^\ast}{\tau_1}\left(i_1^{\ast}
+d_1^{\ast}\right)\left(p_{1,b}^{'}+q_{b}^{'}\right)$\\[3ex]
 $\displaystyle\xi_d=\frac{\partial f_1}{\partial w_{2}(t-\tau_2)}$ &$\displaystyle-\frac{w_1^\ast}{\tau_1}q_{d}^{'}\left(i_1^{\ast}
+d_1^{\ast}\right) $\\[3ex]
 $\displaystyle\xi_{aa}=\frac{1}{2}\frac{\partial^2 f_1}{\partial w_{1}^{2}(t)}$&$\displaystyle\frac{w_1^\ast}{\tau_1}\left(\frac{ i_{1,a}^{''}d_{1}^{*}-i_{1}^{\ast}   d_{1,a}^{''}}{i_1^{\ast}+d_1^{\ast}}\right)$\\[3ex]
 $\displaystyle\xi_{bb}=\frac{1}{2}\frac{\partial^2 f_1}{\partial w_{1}^{2}(t-\tau_1)}$&$\displaystyle-\frac{1}{\tau_1}\left(i_{1}^{\ast}
+d_{1}^{\ast}\right)\left(p_{1,bb}^{''}w_1^\ast+ q_{bb}^{''}w_1^\ast
+ 2p_{1,b}^{'}+2q_{b}^{'}\right)$\\[3ex]
$\displaystyle\xi_{dd}=\frac{1}{2}\frac{\partial^2 f_1}{\partial w_{2}^{2}(t-\tau_2)}$&$\displaystyle  -\frac{w_1^\ast}{\tau_1}q_{dd}^{''}\left(i_{1}^{\ast}
+d_{1}^{\ast}\right)$\\[3ex]
$\displaystyle\xi_{ab}=\frac{\partial^2 f_1}{\partial w_{1}(t)\partial w_{1}(t-\tau_1)}$&$\displaystyle\frac{i_{1,a}^{'}}{\tau_1}-\frac{1}{\tau_1}\left(p_{1}^{\ast}+q^\ast+p_{1,b}^{'}w_1^\ast+q_{b}^{'}w_1^\ast\right)\left(i_{1,a}^{'}+d_{1,a}^{'}
\right)$\\[3ex]
$\displaystyle\xi_{ad}=\frac{\partial^2 f_1}{\partial w_{1}(t)\partial w_{2}(t-\tau_2)}$&$\displaystyle-\frac{w_{1}^{\ast}}{\tau_1}q_{d}^{'}\left(i_{1,a}^{'}+d_{1,a}^{'}
\right)$\\[3ex]
$\displaystyle\xi_{bd}=\frac{\partial^2 f_1}{\partial w_{1}(t-\tau_1)\partial w_{2}(t-\tau_2)}$&$\displaystyle-\frac{1}{\tau_1}\left(i_{1}^{\ast}
+d_{1}^{\ast}\right)\left(q_{d}^{'}+q_{bd}^{''}w_{1}^{\ast}\right)$\\[3ex]
$\displaystyle\xi_{aaa}=\frac{1}{6}\frac{\partial^3 f_1}{\partial w_{1}^{3}(t)}$&$\displaystyle\frac{w_1^\ast}{\tau_1}\left(\frac{ i_{1,a}^{'''}d_{1}^{*}-i_{1}^{\ast}   d_{1,a}^{'''}}{i_1^{\ast}+d_1^{\ast}}\right)$\\[3ex]
$\displaystyle\xi_{bbb}=\frac{1}{6}\frac{\partial^3 f_1}{\partial w_{1}^{3}(t-\tau_1)}$&$\displaystyle-\frac{1}{\tau_1}\left(i_{1}^{\ast}
+d_{1}^{\ast}\right)\left(p_{1,bbb}^{'''}w_1^\ast+ q_{bbb}^{'''}w_1^\ast
+ 3p_{1,bb}^{''}+3q_{bb}^{'}\right)$\\[3ex]
$\displaystyle\xi_{ddd}=\frac{1}{6}\frac{\partial^3 f_1}{\partial w_{2}^{3}(t-\tau_2)}$&$\displaystyle -\frac{w_1^\ast}{\tau_1}q_{ddd}^{'''}\left(i_{1}^{\ast}
+d_{1}^{\ast}\right)$\\[3ex]
$\displaystyle\xi_{aab}=\frac{1}{2}\frac{\partial^3 f_1}{\partial w_{1}^{2}(t)\partial w_{1}(t-\tau_1)}$&$\displaystyle\frac{i_{1,aa}^{''}}{\tau_1}-\frac{1}{\tau_1}\left(p_{1}^{\ast}+q^\ast+p_{1,b}^{'}w_1^\ast+q_{b}^{'}w_1^\ast\right)\left(i_{1
,aa}^{''}+d_{1,aa}^{''}\right)$\\[3ex]
$\displaystyle \xi_{aad}=\frac{1}{2}\frac{\partial^3 f_1}{\partial w_{1}^{2}(t)\partial w_{2}(t-\tau_2)}$&$\displaystyle  -\frac{w_{1}^{\ast}}{\tau_1}q_{d}^{'}\left(i_{1,aa}^{''}+d_{1,aa}^{''}
\right)$\\[3ex]
$\displaystyle\xi_{abb}=\frac{1}{2}\frac{\partial^3 f_1}{\partial w_{1}(t)
\partial w_{1}^{2}(t-\tau_1)}$&$\displaystyle-\frac{1}{\tau_1}\left(i_{1,a}^{'}
+d_{1,a}^{'}\right)\left(p_{1,bb}^{''}w_1^\ast+ q_{bb}^{''}w_1^\ast
+ 2p_{1,b}^{'}+2q_{b}^{'}\right)$\\[3ex]
$\displaystyle \xi_{bbd}=\frac{1}{2}\frac{\partial^3 f_1}{\partial
 w_{1}^{2}(t-\tau_1)\partial w_{2}(t-\tau_2)}$&$\displaystyle-\frac{1}{\tau_1}\left(i_{1}^{\ast}
+d_{1}^{\ast}\right)\left(2q_{bd}^{''}+q_{bbd}^{'''}w_{1}^{\ast}\right)$\\[3ex]
$\displaystyle\xi_{add}=\frac{1}{2}\frac{\partial^3 f_1}{\partial w_{1}(t)\partial w_{2}^{2}(t-\tau_2)}$&$\displaystyle-\frac{w_{1}^{\ast}}{\tau_1}q_{dd}^{''}\left(i_{1,a}^{'}+d_{1,a}^{'}
\right) $\\[3ex]
$\displaystyle\xi_{bdd}=\frac{1}{2}\frac{\partial^3 f_1}{\partial w_{1}(t-\tau_1)\partial w_{2}^{2}(t-\tau_2)}$&$\displaystyle-\frac{1}{\tau_1}\left(i_{1}^{\ast}
+d_{1}^{\ast}\right)\left(q_{dd}^{''}+q_{bdd}^{'''}w_{1}^{\ast}\right)$\\[3ex]
$\displaystyle\xi_{abd}=\frac{\partial^3 f_1}{\partial w_{1}(t)\partial w_{1}(t-\tau_1)\partial w_{2}(t-\tau_2)}$&$\displaystyle-\frac{1}{\tau_1}\left(i_{1,a}^{'}
+d_{1,a}^{'}\right)\left(q_{bd}^{''}w_1^\ast
+q_{d}^{'}\right)$\\[3ex]
 $\displaystyle\chi_c=\frac{\partial f_2}{\partial w_{2}(t)}$ & $\displaystyle\frac{w_2^\ast}{\tau_2}\left(\frac{ i_{2,c}^{'}d_{2}^{*}-i_{2}^{\ast}   d_{2,c}^{'}}{i_2^{\ast}+d_2^{\ast}}\right)$ \\[3ex]
\hline\hline
\label{table:coefficients}
\end{tabular}
\end{table*}

\begin{table*}
\centering
\setlength{\tabcolsep}{36pt}
\begin{tabular}{l l }\hline\hline
\\ [-1mm]
 $\displaystyle\chi_d=\frac{\partial f_2}{\partial w_{2}(t-\tau_2)}$ & $\displaystyle -\frac{w_2^\ast}{\tau_2}\left(i_2^{\ast}
+d_2^{\ast}\right)\left(p_{2,d}^{'}+q_{d}^{'}\right)$\\[3ex]
$\displaystyle\chi_b=\frac{\partial f_2}{\partial w_{1}(t-\tau_1)}$ &$\displaystyle-\frac{w_2^\ast}{\tau_2}q_{b}^{'}\left(i_2^{\ast}
+d_2^{\ast}\right) $\\[3ex]
 $\displaystyle\chi_{cc}=\frac{1}{2}\frac{\partial^2 f_2}{\partial w_{2}^{2}(t)}$&$\displaystyle\frac{w_2^\ast}{\tau_2}\left(\frac{ i_{2,c}^{''}d_{2}^{*}-i_{2}^{\ast}   d_{2,c}^{''}}{i_2^{\ast}+d_2^{\ast}}\right)$\\[3ex]
 $\displaystyle\chi_{dd}=\frac{1}{2}\frac{\partial^2 f_2}{\partial w_{2}^{2}(t-\tau_2)}$&$\displaystyle-\frac{1}{\tau_2}\left(i_{2}^{\ast}
+d_{2}^{\ast}\right)\left(p_{2,dd}^{''}w_2^\ast+ q_{dd}^{''}w_2^\ast
+ 2p_{2,d}^{'}+2q_{d}^{'}\right)$\\[3ex]
$\displaystyle\chi_{bb}=\frac{1}{2}\frac{\partial^2 f_2}{\partial w_{1}^{2}(t-\tau_1)}$&$\displaystyle  -\frac{w_2^\ast}{\tau_2}q_{bb}^{''}\left(i_{2}^{\ast}
+d_{2}^{\ast}\right)$\\[3ex]
$\displaystyle\chi_{cd}=\frac{\partial^2 f_2}{\partial w_{2}(t)\partial w_{2}(t-\tau_2)}$&$\displaystyle\frac{i_{2,c}^{'}}{\tau_2}-\frac{1}{\tau_1}\left(p_{2}^{\ast}+q^\ast+p_{2,d}^{'}w_2^\ast+q_{d}^{'}w_2^\ast\right)
\left(i_{2,c}^{'}+d_{2,c}^{'}
\right)$\\[3ex]
$\displaystyle\chi_{bc}=\frac{\partial^2 f_2}{\partial w_{1}(t-\tau_1)\partial w_{2}(t)}$&$\displaystyle-\frac{w_{2}^{\ast}}{\tau_2}q_{b}^{'}\left(i_{2,c}^{'}+d_{2,c}^{'}
\right)$\\[3ex]
$\displaystyle\chi_{bd}=\frac{\partial^2 f_2}{\partial w_{1}(t-\tau_1)\partial w_{2}(t-\tau_2)}$&$\displaystyle-\frac{1}{\tau_2}\left(i_{2}^{\ast}
+d_{2}^{\ast}\right)\left(q_{b}^{'}+q_{bd}^{''}w_{2}^{\ast}\right)$\\[3ex]
$\displaystyle\chi_{ccc}=\frac{1}{6}\frac{\partial^3 f_2}{\partial w_{2}^{3}(t)}$&$\displaystyle\frac{w_2^\ast}{\tau_2}\left(\frac{ i_{2,c}^{'''}d_{2}^{*}-i_{2}^{\ast}   d_{2,c}^{'''}}{i_2^{\ast}+d_2^{\ast}}\right)$\\[3ex]
$\displaystyle\chi_{ddd}=\frac{1}{6}\frac{\partial^3 f_2}{\partial w_{2}^{3}(t-\tau_2)}$&$\displaystyle-\frac{1}{\tau_2}\left(i_{2}^{\ast}
+d_{2}^{\ast}\right)\left(p_{2,ddd}^{'''}w_2^\ast+ q_{ddd}^{'''}w_2^\ast
+ 3p_{2,dd}^{''}+3q_{dd}^{'}\right)$\\[3ex]
$\displaystyle\chi_{bbb}=\frac{1}{6}\frac{\partial^3 f_2}{\partial w_{1}^{3}(t-\tau_1)}$&$\displaystyle -\frac{w_2^\ast}{\tau_2}q_{bbb}^{'''}\left(i_{2}^{\ast}
+d_{2}^{\ast}\right)$\\[3ex]
$\displaystyle\chi_{ccd}=\frac{1}{2}\frac{\partial^3 f_2}{\partial w_{2}^{2}(t)\partial w_{2}(t-\tau_2)}$&$\displaystyle\frac{i_{2,cc}^{''}}{\tau_2}-\frac{1}{\tau_2}\left(p_{2}^{\ast}+q^\ast+p_{2,d}^{'}w_2^\ast+q_{d}^{'}w_2^\ast\right)
\left(i_{2,cc}^{''}+d_{2,cc}^{''}\right)$\\[3ex]
$\displaystyle \chi_{bcc}=\frac{1}{2}\frac{\partial^3 f_2}{\partial w_{1}(t-\tau_1)\partial w_{2}^{2}(t)}$&$\displaystyle  -\frac{w_{2}^{\ast}}{\tau_2}q_{b}^{'}\left(i_{2,cc}^{''}+d_{2,cc}^{''}
\right)$\\[3ex]
$\displaystyle\chi_{cdd}=\frac{1}{2}\frac{\partial^3 f_2}{\partial w_{2}(t)
\partial w_{2}^{2}(t-\tau_2)}$&$\displaystyle-\frac{1}{\tau_2}\left(i_{2,c}^{'}
+d_{2,c}^{'}\right)\left(p_{2,dd}^{''}w_2^\ast+ q_{dd}^{''}w_2^\ast
+ 2p_{2,d}^{'}+2q_{d}^{'}\right)$\\[3ex]
$\displaystyle \chi_{bdd}=\frac{1}{2}\frac{\partial^3 f_2}{\partial
 w_{1}(t-\tau_1)\partial w_{2}^{2}(t-\tau_2)}$&$\displaystyle-\frac{1}{\tau_2}\left(i_{2}^{\ast}
+d_{2}^{\ast}\right)\left(2q_{bd}^{''}+q_{bdd}^{'''}w_{2}^{\ast}\right)$\\[3ex]
$\displaystyle\chi_{bbc}=\frac{1}{2}\frac{\partial^3 f_2}{\partial w_{1}^{2}(t-\tau_1)\partial w_{2}(t)}$&$\displaystyle-\frac{w_{2}^{\ast}}{\tau_2}q_{bb}^{''}\left(i_{2,c}^{'}+d_{2,c}^{'}
\right) $\\[3ex]
$\displaystyle\chi_{bbd}=\frac{1}{2}\frac{\partial^3 f_1}{\partial w_{1}^{2}(t-\tau_1)\partial w_{2}(t-\tau_2)}$&$\displaystyle-\frac{1}{\tau_2}\left(i_{2}^{\ast}
+d_{2}^{\ast}\right)\left(q_{bb}^{''}+q_{bbd}^{'''}w_{2}^{\ast}\right)$\\[3ex]
$\displaystyle\xi_{bcd}=\frac{\partial^3 f_1}{\partial w_{1}(t-\tau_1)\partial w_{2}(t)\partial w_{2}(t-\tau_2)}$&$\displaystyle-\frac{1}{\tau_2}\left(i_{2,c}^{'}
+d_{2,c}^{'}\right)\left(q_{bd}^{''}w_2^\ast
+q_{b}^{'}\right)$\\[3ex]
 \hline\hline
 \end{tabular}
\end{table*}

\begin{thebibliography}{100}


\bibitem{Cerf}
V.G.~Cerf, ``Bufferbloat and other Internet challenges", \emph{IEEE Internet Computing}, vol. 5, pp. 79--80, 2014.

\bibitem{Cooke}
K.L.~Cooke, and Z.~Grossman, ``Discrete delay, distributed delay and stability switches", \emph{Journal of Mathematical Analysis and Applications}, vol.~86, pp.~592--627, 1982.

\bibitem{DDE1}
K.~Engelborghs, T.~Luzyanina, and D.~Roose, ``Numerical bifurcation analysis of delay differential equations using DDE-BIFTOOL", \emph{ACM Transactions on Mathematical Software (TOMS}), vol.~28, pp.~1--21, 2002.


\bibitem{DDE2}
K.~Engelborghs, T.~Luzyanina, G.~Samaey, ``DDE-BIFTOOL v. 2.00: a Matlab package for bifurcation analysis of delay differential equations", \emph{Technical Report TW-330, Department of Computer Science, K.U.Leuven, Leuven, Belgium}, 2001. 


\bibitem{Gettys}
J.~Gettys and K.~Nichols, ``Bufferbloat: dark buffers in the
Internet", \emph{Communications of the ACM}, vol. 55, pp. 57--65,
2012.
\bibitem{Jagannathan}
D.~Ghosh, K.~Jagannathan, and G.~Raina, ``Right buffer sizing matters: stability, queuing delay and traffic burstiness in compound TCP", in \emph{Proceedings of 52nd Annual Allerton Conference on Communication, Control, and Computing}, 2014.


\bibitem{Ha}
S.~Ha, I.~Rhee and L.~Xu, ``CUBIC: a new TCP-friendly high-speed TCP variant",
\emph{ACM SIGOPS Operating Systems Review}, vol. 42, pp. 64--74, 2008.


\bibitem{Hassard}
B.D.~Hassard, N.D.~Kazarinoff and Y-H.~Wan, \emph{Theory
and Applications of Hopf Bifurcation}, Cambridge University
Press, 1981.
\bibitem{Kuznetsov}
Y.A.~Kuznetsov, \emph{Elements of Applied Bifurcation Theory}, Springer Science \& Business Media, 2013.

\bibitem{Nichols}
K.~Nichols and V.~Jacobson, ``Controlling queue delay",
\emph{Communications of the ACM}, vol. 55, pp. 42--50, 2012.

\bibitem{Padhye}
J.~Padhye, V.~Firoiu, D.~Towsley and J.F.~Kurose, ``Modeling TCP Reno performance: a simple model and its empirical validation", \emph{IEEE/ACM Transactions
on Networking}, vol. 8, pp. 133--145, 2000.


\bibitem{Gaurav}
G. Raina, ``Local bifurcation analysis of some dual congestion control algorithms", \emph{IEEE Transactions on Automatic
Control}, vol. 50, pp. 1135--1146, 2005.

\bibitem{Raja}
P. Raja and G. Raina, ``Delay and loss-based transport protocols: buffer-sizing and stability", in \emph{Proceedings of International Conference on Communication Systems and Networks}, 2012.
\bibitem{Rudin}
W.~Rudin, \emph{Real and Complex Analysis}, Tata McGraw-Hill Education, 1987.


\bibitem{Tan}
K.~Tan, J.~Song, Q.~Zhang and M.~Sridharan, ``A Compound
TCP approach for high-speed and long distance networks",
in \emph{Proceedings of IEEE INFOCOM}, 2006.


\bibitem{ns2}
The Network Simulator (NS2). [Online]. Available:
http://nsnam.isi.edu/nsnam/index.php/User Information.





\end{thebibliography}
\end{document}